\documentclass{amsproc}
\usepackage{amsmath,amsfonts,amsthm,mathtools}
\usepackage{geometry}

\usepackage{mathtools}
\usepackage{enumitem}
\usepackage{xcolor}
\usepackage{tikz}
\usepackage{graphicx}
\usepackage{epstopdf}
\newcommand{\R}{\mathbb{R}}
\newcommand{\ve}{\varepsilon}
\newcommand{\cO}{\mathcal{O}}
\newcommand{\cC}{\mathcal{C}}
\newcommand{\cS}{\mathcal{S}}
\newcommand{\cB}{\mathcal{B}}
\newcommand{\cM}{\mathcal{M}}

\newcommand{\bbS}{\mathbb{S}}
\newcommand{\bbN}{\mathbb{N}}

\newcommand{\txta}{\textnormal{a}}
\newcommand{\txtr}{\textnormal{r}}
\newcommand{\EN}{\textnormal{en}}
\newcommand{\EX}{\textnormal{ex}}
\newcommand{\ddt}[1]{\frac{\textnormal{d}#1}{\textnormal{d}t}}

\newcommand{\D}{\textnormal{D}}

\newtheorem{theorem}{Theorem}[section]
\newtheorem{lemma}[theorem]{Lemma}

\theoremstyle{definition}
\newtheorem{definition}[theorem]{Definition}

\theoremstyle{remark}
\newtheorem{remark}[theorem]{Remark}

\numberwithin{equation}{section}

\newcommand{\br}{\bar{r}}
\newcommand{\bx}{\bar{x}}
\newcommand{\by}{\bar{y}}
\newcommand{\be}{\bar{\ve}}
\newcommand{\ch}[1]{#1}
\newcommand{\chh}[1]{#1}
  
\usepackage{hyperref} 
\begin{document}
\title{A survey on the blow-up method for fast-slow systems}

\author{Hildeberto Jard\'on-Kojakhmetov}
\address{Technical University of Munich, Faculty of Mathematics, Research Unit ``Multiscale and Stochastic Dynamics'', 85748 Garching b. M{\"u}nchen, Germany.}
\email{h.jardon.kojakhmetov@tum.de}
\thanks{$^*$~Invited survey article.}

\author{Christian Kuehn}
\address{Technical University of Munich, Faculty of Mathematics, Research Unit ``Multiscale and Stochastic Dynamics'', 85748 Garching b. M{\"u}nchen, Germany.}
\email{ckuehn@ma.tum.de}
\thanks{}

\begin{abstract}
	In this document we review a geometric technique, called \emph{the blow-up method}, as it has been used to analyze and understand the dynamics of fast-slow systems around non-hyperbolic points. The blow-up method, having its origins in algebraic geometry, was introduced in 1996 to the study of fast-slow systems in the seminal work by Dumortier and Roussarie \cite{dumortier1996canard}, whose aim was to give a geometric approach and interpretation of canards in the van der Pol oscillator. Following \cite{dumortier1996canard}, many efforts have been performed to expand the capabilities of the method and to use it in a wide range of scenarios. Our goal is to present in a concise and compact form those results that, based on the blow-up method, are now the foundation of the geometric theory of fast-slow systems with non-hyperbolic singularities. Due to their great importance in the theory of fast-slow systems, we cover fold points as one of the main topics. Furthermore, we also present several other singularities such as Hopf, pitchfork, transcritical, cusp, and Bogdanov-Takens, in which the blow-up method has been proved to be extremely useful. Finally, we survey further directions as well as examples of specific applied models, where the blow-up method has been used successfully.
\end{abstract}

\maketitle

\tableofcontents

\section{Introduction}

A large number of phenomena in nature can be understood as the interaction of  simpler sub-processes. There are cases where it can be observed that such sub-processes evolve at distinct time scales. A classical example is found in enzymatic reactions (chemical reactions that are accelerated or catalysed through the addition of a particular type of protein called enzyme), where it can be often measured, and thus assumed in the mathematical models, that the enzyme's concentration evolves and reaches its final state much faster than any of the other reactants \cite{briggs1925note}. Other examples of phenomena where different time scales evolutions may be observed can be found  in neuron dynamics, population dynamics, network dynamics (such as social and power networks), biological processes (such as cell division and regulatory processes)~\cite{banasiak2014methods,hek2010geometric}, among many others. Mathematical models in the form of Ordinary Differential Equations (ODEs) that reflect the presence of multiple time scales are called \emph{fast-slow systems} or \emph{singularly perturbed ODEs}. 

From a mathematical perspective, fast-slow systems have been studied for a long time and several techniques including asymptotic and geometric methods~\cite{eckhaus2011matched,kevorkian2012multiple,kuehn2015multiple,OMalley1991,verhulst2005methods} have been developed to understand them.

In qualitative terms, the key idea when studying fast-slow systems is to separate the sub-processes acting at the different time scales, understand them, and then try to describe the full dynamics based on the subsystems. \ch{The previous rough idea can be made mathematically rigorous and, in fact, is the basis of the geometric analysis of fast-slow systems as we describe in Section \ref{sec:preliminaries}.} \ch{However, as we shall detail in the main part of this paper, there are cases in which the fast and slow dynamics cannot be distinguished, and thus separated, in a clear and straightforward manner. At the level of the mathematical model, this obstacle carries several difficulties, the most common is the presence of non-hyperbolic equilibrium points (or singularities), and thus advanced methods of analysis are required in such a case. The presence or absence of singularities determines what type of mathematical technique is suitable for analysis. This article is dedicated to review a technique, called \emph{the blow-up method}, that allows one to study fast-slow systems near \chh{nilpotent} singularities.}

\chh{Briefly speaking, the blow-up method is a change of coordinates, see the details in Sections \ref{sec:bu} and \ref{sec:blowup}, with which one transforms the analysis of fast-slow systems around nilpotent singularities to a series of sub-problems with (semi-)hyperbolic singularities, which can be dealt with known methods of dynamical systems}. Afterwards, we gather all the information obtained from the sub-problems and draw conclusions on the dynamics of the full fast-slow system. It is important to remark that the blow-up (or also blowing up) method was initially developed for the resolution of singularities of algebraic curves, and more generally of varieties over fields of characteristic zero~\cite{Hironaka1,Hironaka2}. The idea behind resolution of singularities consists on a certain \ch{coordinate mapping} that provides a regular curve from a singular one. We will see that this idea is similar for ordinary differential equations. In the context of dynamical systems, the blow-up method has been used, for example, to study and classify \ch{nilpotent} singularities of (single-scale) vector fields~\cite{alvarez2011survey,brunella1990topological,dumortier1977singularities,dumortier1991local,takens1974singularities}, see also Section \ref{sec:bu} below. Here we focus on the use of blow-up for fast-slow multiscale systems.

Our goal with this survey is to gather in a concise and compact form the most relevant results in the theory of fast-slow systems for which the blow-up method has been fundamental. The results that we include in the main part of this survey have extended and complemented the general theory of fast-slow systems and are nowadays the foundation of many other studies. Our approach is structured as follows: in Section \ref{sec:preliminaries} we provide some necessary mathematical preliminaries, where we sketch Fenichel's theorem and the blow-up method. Later, in Section \ref{sec:survey} we present our survey which includes results on generic folds and canards. Then we recall some fundamental results outside the scope of folded singularities, and present several applications, where the blow-up method has been used. We conclude in Section \ref{sec:summary} with a summary and an outlook for the blow-up method within the context of fast-slow systems.

\begin{center}
	\vspace*{.25cm}
	\textbf{Notation}
\end{center}

Naturally, different notations are used across the whole body of literature regarding fast-slow systems. In this article we homogenize the notation used for the most important geometric objects that one encounters while studying fast-slow systems. Although the explicit definition of the geometric objects changes depending on the specific problem, their role in the analysis is often similar. Clarifying the way we decide to denote these objects allows us to recycle the notation, thus facilitating the exposition. The reader may refer back to the following conventions throughout:

\begin{itemize}[leftmargin=*]
	\item $\cC_0$ denotes the critical manifold of a given fast-slow system to be specified when appropriate.
	\item $\cS_0$ is a compact subset of $\cC_0$. The attracting and repelling parts of $\cS_0$ are denoted by $\cS_0^\txta$ and $\cS_0^\txtr$ respectively. When further distinction is needed, we add appropriate subscripts.
	\item  $\cS_\ve$ denotes a slow manifold, and $\cS_\ve^\txta$ (resp. $\cS_\ve^\txtr$) denotes the attracting (resp. repelling) part of $\cS_\ve$.
	\item $\Sigma^\EN$ is a codimension $1$ \ch{(that is, the dimension of $\Sigma^\EN$ is one less as that of the phase-space)} \emph{entry} section of the phase-space transverse to $\cS_0^\txta$. $\Sigma^\EX$ is used to denote a codimension $1$ \emph{exit} section of the phase-space transverse to $\cS_\ve^\txta$ beyond a non-hyperbolic singularity. When more than one exit section is necessary, we add appropriate subscripts.
	\item $\Pi$ denotes the transition map between $\Sigma^\EN$ and $\Sigma^\EX$ induced by the flow of a fast-slow system. When more than one transition map is necessary, we add appropriate subscripts.
	\item $K_\bullet$ denotes a chart of the blow-up space. We distinguish the charts using a subscript depending on the number of charts required for the analysis of a specific singularity.
\end{itemize}

\section{Preliminaries}\label{sec:preliminaries}

\ch{
	In this section we provide the mathematical preliminaries needed for the survey given in Section \ref{sec:survey}. First we describe the main idea of the blow-up method for a single time scale planar vector field. Next, we state a couple of important definitions regarding fast-slow systems. Afterwards, we recall the fundamental Fenichel's theorem, and finish this section with a description of the blow-up method for fast-slow systems.
}
	
	\subsection{\ch{The blow-up method}} \label{sec:bu}
	\ch{
	In this section we introduce the blow-up method in its classical context, that is, to desingularize a nilpotent equilibrium point\footnote{We recall that an equilibrium point of a vector field is called nilpotent if the linearization of the vector field at such a point is given by a matrix with only zero eigenvalues.} of a planar vector field. For a detailed exposition the reader is referred to \cite{alvarez2011survey,dumortier1991local} and \cite[Chapter 7]{kuehn2015multiple}. Here we shall only treat an example to highlight the main idea of the method. Later, in Section \ref{sec:blowup} we will see how this transformation also fits into the study of fast-slow systems.

	Let us consider the planar ordinary differential equation (ODE)
	\begin{equation}\label{eq:ex1}
		\begin{split}
			\ddt{x} &= y\\
			\ddt{y} &= x^3 + xy.
		\end{split}
	\end{equation}
	We note that the origin $(x,y)=(0,0)$ is a unique equilibrium point and that the linearization of \eqref{eq:ex1} at the origin is given by the matrix
	\begin{equation}
		\begin{bmatrix}
			0 & 1\\ 0 & 0
		\end{bmatrix}.
	\end{equation}
	
	Thus, the origin is a non-hyperbolic equilibrium point and, moreover, is nilpotent. Our goal is to qualitatively describe the orbits of \eqref{eq:ex1} in a small neighbourhood of the origin. However, not only the linearization offers no useful information, but centre manifold reduction \cite{carr2012applications} is not suitable since in this case the centre manifold corresponds to the whole phase-space. \chh{So, what we are going to use is a suitable change of coordinates, known as \emph{blow-up}, which will induce a new system with only hyperbolic equilibrium points, and therefore can be analyzed by dynamical systems tools}.

	Let us consider a \emph{weighted polar change of coordinates}
	\begin{equation}\label{eq:ex-bu}
		\phi:\bbS^1\times I \to\R^2, \qquad \qquad \phi(\theta,r)=(r\cos\theta,r^2\sin\theta),
	\end{equation}
	where $I\subseteq\R$ is an interval containing the origin and $\theta\in[0,2\pi]$. At the end of this section we clarify the reason to choose a weighted polar change of coordinates, for now let us proceed with the example.

	The change of coordinates defined by $(x,y)=(r\cos\theta,r^2\sin\theta)$ defines a new ODE, namely
	\begin{equation}\label{eq:ex2}
		\begin{split}
			\dot\theta &= \frac{r}{\sin^2\theta+1}\left(1+\sin\theta-4\sin^2\theta-\sin^3\theta+\sin^4\theta \right)\\
			\dot r &= \frac{r^2}{\sin^2\theta+1}\cos\theta\sin\theta\left( \sin\theta-\sin^2\theta+2 \right).
		\end{split}
	\end{equation}

	Note that the change of coordinates defined by $\phi$ maps the circle $\bbS^1\times\left\{ 0 \right\}$ to the origin in the plane\footnote{Equivalently $\phi^{-1}$ maps the origin in the plane to the circle $\bbS^1\times\left\{ 0 \right\}$.}. Moreover, since $\phi$ is a diffeomorphism for $\left\{ r>0 \right\}$,  orbits of \eqref{eq:ex1} in a small neighbourhood of the origin correspond to orbits of \eqref{eq:ex2} in a small neighbourhood of $\bbS^1\times\left\{ 0 \right\}$. Note however that \eqref{eq:ex2} vanishes along $\bbS^1\times\left\{ 0 \right\}$. To overcome this we can divide the right-hand side of \eqref{eq:ex2} by $r$. This operation does not change the qualitative properties of the orbits in the region $\bbS^1\times\left\{ r>0 \right\}$. Thus, it shall suffice to study \emph{the desingularized system}
	\begin{equation}\label{eq:ex3}
		\begin{split}
			\dot\theta &= \frac{1}{\sin^2\theta+1}\left(1+\sin\theta-4\sin^2\theta-\sin^3\theta+\sin^4\theta \right)\\
			\dot r &= \frac{r}{\sin^2\theta+1}\cos\theta\sin\theta\left( \sin\theta-\sin^2\theta+2 \right),
		\end{split}
	\end{equation}
	which does not vanish any more along $\bbS^1\times\left\{ 0 \right\}$. The most important fact is that \emph{orbits of \eqref{eq:ex3} near $\bbS^1\times\left\{ 0 \right\}$ correspond to orbits of \chh{\eqref{eq:ex1}} near the origin.}

	It is now straightforward to show that \eqref{eq:ex3} has four \emph{hyperbolic} saddle equilibrium points, namely $p_1=(-\arcsin(\sqrt{2}-1),0)$, $p_2=(\arcsin(\sqrt{5}/2-1/2),0)$, $p_3=(\pi-\arcsin(\sqrt{5}/2-1/2),0)$ and $p_4=(\pi+\arcsin(\sqrt{2}-1),0)$. Since the aforementioned equilibrium points are hyperbolic it follows from linear analysis that the phase portrait of \eqref{eq:ex3} in a small neighbourhood of $\bbS^1\times\left\{ 0 \right\}$ is as show in Figure \ref{fig:ex1}.

	\begin{figure}[htbp]\centering
	\begin{tikzpicture}
		\node at (-6,0){
		\includegraphics[scale=1.25]{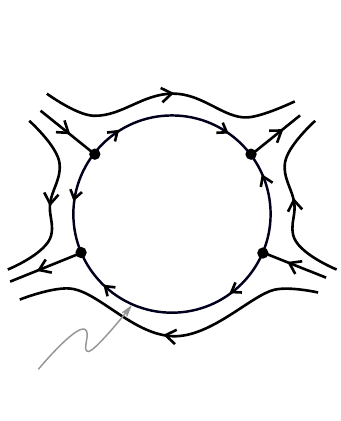}
		};
		\node at (3,0){
		\includegraphics[scale=1.25]{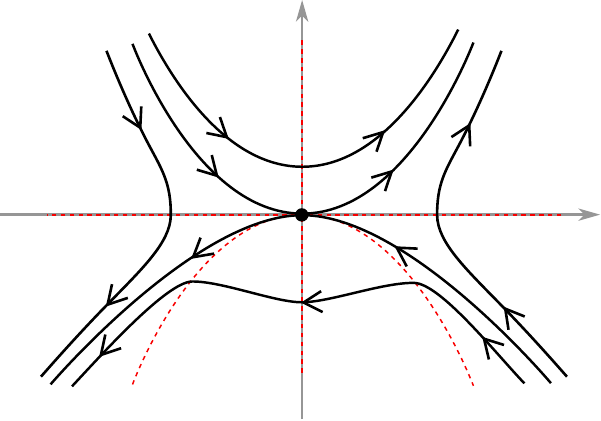}
		};
		\node at (-7.5,-2.3){$\bbS^1\times\left\{ 0 \right\}$};
		\node at (-5.1,-.5){$p_1$};
		\node at (-5.1, .5){$p_2$};
		\node at (-6.7, .5){$p_3$};
		\node at (-6.8,-.5){$p_4$};

		\node at (6.9,-.075){$x$};
		\node at (3,2.9){$y$};

		\draw[black,thick,<-] (-1.5,2) arc (60:120:2) node[midway,above]{$\phi$};
		\draw[black,thick,->] (-1.5,-2) arc (-60:-120:2) node[midway,below]{$\phi^{-1}$};
	\end{tikzpicture}
		  \label{fig:ex1}
		  \caption{Blow-up analysis of \eqref{eq:ex1}. On the left we show the phase-portrait of \eqref{eq:ex3} in a small neighborhood of $\bbS^1\times\left\{ 0 \right\}$, where four hyperbolic saddle points are found. On the right we show the corresponding orbits of \eqref{eq:ex1}, where from a qualitative perspective, the circle $\bbS^1\times\left\{ 0 \right\}$ ``blows-down'' to the origin and all other orbits of \eqref{eq:ex1} are equivalent to orbits of \eqref{eq:ex3}. To provide more detail on the flow of \eqref{eq:ex1} away from the origin we have made use of the corresponding nullclines, shown as dashed-red curves.}
	\end{figure}
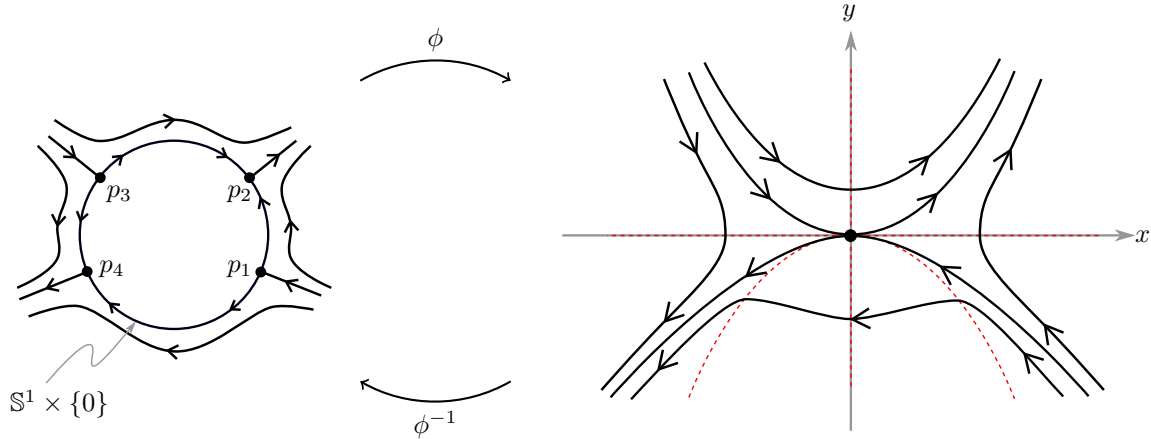

We finish this section with some important remarks:
\begin{itemize}[leftmargin=*]
	\item The procedure we exemplified above is known as \emph{the blow-up method}. In some sense, the transformation $\phi^{-1}$ ``blows the origin up to a circle''. The advantage of blowing up is that one obtains a new system which is simpler to analyze. We recall that, in the above example, \eqref{eq:ex1} has a nilpotent equilibrium point at the origin while \eqref{eq:ex3} has four hyperbolic equilibrium points along $\bbS^1\times\left\{ 0 \right\}$, which are simpler to study with standard techniques of dynamical systems. Once the blown-up system is understood we then ``blow-down'' the phase-portrait of \eqref{eq:ex3} resulting in a qualitative description of the original system \eqref{eq:ex1}. 
	\item In the example presented above we have used a weighted version of a polar change of coordinates. Usually one then refers to the transformation as a \emph{quasi-homogeneous} blow-up to emphasize that the weights in the transformation are distinct from $1$. The advantage of using a quasi-homogeneous blow-up instead of a homogeneous one is that we can desingularize the origin in just one step. The reader can check that if one uses $(x,y)=(r\cos\theta,r\sin\theta)$ instead of \eqref{eq:ex-bu}, the blown-up system then has a pair of nilpotent singularities located at $(\theta,r)=(0,0)$ and $(\theta,r)=(\pi,0)$. In turn, the blow-up method can be applied once more to such pair of points. More details are provided in Section \ref{sec:blowup}, see also \cite{dumortier1991local} and \cite[Chapter 7]{kuehn2015multiple}.
	\item It is known, for example, that real-analytic planar and three dimensional vector fields can be desingularized via the blow-up method \cite{dumortier1977singularities,panazzolo2002desingularization,panazzolo2006resolution}. In particular, several of the fast-slow systems in this survey fall into this category. The blow-up method has been used in higher dimensional fast-slow systems only in a case-by-case approach. 
\end{itemize}

	}

	\subsection{\ch{Fast-slow systems}}

A fast-slow system is a singularly perturbed ordinary differential equation (ODE) of the form
\begin{equation}\label{eq:sfs1}
	\begin{split}
		\ve\dot x &= f(x,y,\ve),\\
		\dot y &= g(x,y,\ve),
	\end{split}
\end{equation}
where the over-dot denotes the derivative with respect to \emph{the slow time} $\tau$, $x\in\R^m$ and $y\in\R^n$ denote the fast and the slow variables respectively, and $\ve>0$ is a small parameter accounting for the time scale separation. \ch{System \eqref{eq:sfs1} is said to be \emph{in standard form}, and we can explicitly see the time scale separation between the slow and fast variables. In some cases, however, fast-slow systems may appear in a more general non-standard form, namely
\begin{equation}
	\begin{split}
		\dot\zeta=F(\zeta,\ve),
	\end{split}
\end{equation}
where the time scale separation, if any, may not be evident to distinguish, see Section \ref{sec:applications} for some examples of fast-slow systems which are not in standard form. In this survey we shall adhere to \eqref{eq:sfs1}.
}

By defining the fast time $t\coloneqq \tau/\ve$, we can rewrite \eqref{eq:sfs1} as an $\ve$-family of ODEs of the form
\begin{equation}\label{eq:sfs2}
	\begin{split}
		x' &= f(x,y,\ve),\\
		y' &= \ve g(x,y,\ve),
	\end{split}
\end{equation}
where now the prime $'$ denotes the derivative with respect to \emph{the fast time} $t$. We assume that \ch{$f:\R^m\times\R^n\times\R\to\R^m$ and $g:\R^m\times\R^n\times\R\to\R^n$ are of class $C^k$}, with $k$ sufficiently large. We observe that \eqref{eq:sfs1} and \eqref{eq:sfs2} are smoothly equivalent for $\ve>0$, meaning that their only difference is the time parametrization of the corresponding trajectories.

The term Geometric Singular Perturbation Theory (GSPT)~\cite{fenichel1979geometric,jones1995geometric}, refers to  the collection of geometric techniques with which fast-slow systems can be analyzed. Roughly speaking, the idea of GSPT is to study the limits of \eqref{eq:sfs1} and \eqref{eq:sfs2} as $\ve\to0$, analyze invariant objects in the two limits, and then use perturbation methods to describe the dynamics of \eqref{eq:sfs1} and \eqref{eq:sfs2} for $\ve>0$ sufficiently small. Accordingly, by setting $\ve=0$ in \eqref{eq:sfs1} we obtain
\begin{equation}\label{eq:cde}
	\begin{split}
		0 &= f(x,y,0),\\
		\dot y &= g(x,y,0),
	\end{split}
\end{equation}
which is called the \emph{reduced system} (or \emph{slow subsystem}), and the flow of \eqref{eq:cde} is called \emph{the slow flow}. Structurally, the system~\eqref{eq:cde} is a \emph{Constrained Differential Equation} (CDE)~\cite{takens1976constrained} (or \emph{Differential Algebraic Equation}~\cite{kunkel2006differential}). From~\eqref{eq:sfs1} we get for $\ve=0$ the ODE
\begin{equation}\label{eq:layer}
	\begin{split}
		x' &= f(x,y,0),\\
		y' &=0,
	\end{split}
\end{equation}
which is called \emph{the layer equation} (or \emph{fast subsystem}) and where we can view the slow variables $y$ as parameters. The two \emph{singular limit} systems \eqref{eq:cde} and \eqref{eq:layer} are not equivalent any more. However, we observe that the set $\left\{ f(x,y,0)=0 \right\}$ defines the phase-space of solutions of the slow subsystem and the set of equilibrium points of the fast subsystem equation. This leads to the following natural definition:

\begin{definition}[Critical manifold] The critical manifold is defined as the set
\begin{equation}
	\cC_0 = \left\{ (x,y)\in\R^m\times\R^n\, | \, f(x,y,0)=0 \right\}.
\end{equation}
\end{definition}

\ch{Although we call $\cC_0$ a manifold, strictly speaking, solutions of the equation $f(x,y,0)=0$ do not necessarily define a manifold, that is, they may not be locally diffeomorphic to an Euclidean space. For instance, the solutions of $f(x,y,0)=0$ may self-intersect as in the case of the transcritical and pitchfork singularities in Section \ref{sec:beyondfold} }. However, generically, there are open regions in the phase-space where $\cC_0$ is indeed a manifold. So, as a convention in fast-slow systems one keeps referring to $\cC_0$ as the critical manifold. An important characteristic that critical manifolds may posses is normal hyperbolicity:

\begin{definition}\label{def:nh}
	A point $p=(x,y)\in \cC_0$ is called \emph{hyperbolic} if the matrix $\D_xf(p,0)\in\R^{m\times m}$, where $\D_x$ stands for the total derivative with respect to $x$,  has no eigenvalues on the imaginary axis. The critical manifold $\cC_0$ is called \emph{normally hyperbolic} if all points $p\in\cC_0$ are hyperbolic.
\end{definition}

Note that equivalently to Definition \ref{def:nh}, one can say that $\cC_0$ is normally hyperbolic if all its points are hyperbolic equilibrium points of the fast subsystem. We call a hyperbolic point $p$ \emph{attracting} if all eigenvalues have negative real part, \emph{repelling} if all eigenvalues have positive real part, and \emph{saddle-type} if there are negative and positive real part eigenvalues. On the contrary situation, we shall call $p\in \cC_0$ \emph{non-hyperbolic} if the matrix $\D_x(p,0)$ has at least one eigenvalue on the imaginary axis. Whether the critical manifold is, or is not, normally hyperbolic will distinguish the tools that are needed for the analysis of the corresponding fast-slow system. In the case of normal hyperbolicity we have Fenichel's theorem (see Section \ref{sec:fenichel}); in the case of loss of normal hyperbolicity due to the presence of non-hyperbolic points, we may use the blow-up method (see Section \ref{sec:blowup}). \ch{We emphasize that loss of normal hyperbolicity in fast-slow systems is highly relevant as it is related to dynamic features such as relaxation oscillations, canards, bursting, mixed-mode oscillations, etc\'etera \cite{desroches2012mixed,kuehn2015multiple}.}

\subsection{Fenichel's theorem}\label{sec:fenichel}

Fenichel's theorem~\cite{fenichel1979geometric} (see also~\cite{jones1995geometric,tikhonov1952systems}) provides geometric tools and techniques to analyze fast-slow systems with hyperbolic points:

\begin{theorem}\label{thm:fenichel}
	Suppose that $\cS_0\subseteq\cC_0$ is a compact normally hyperbolic submanifold (possibly with boundary) of the critical manifold $\cC_0$ of \eqref{eq:sfs2}. Then, for $\ve>0$ sufficiently small, the following hold:
	\begin{enumerate}[leftmargin=*]
		\item[(M1)] There exists a locally invariant manifold $\cS_\ve$ diffeomorphic to $\cS_0$. Local invariance means that trajectories can enter or leave $\cS_\ve$ only through its boundaries.
		\item[(M2)] $\cS_\ve$ has Hausdorff distance $\cO(\ve)$, as $\ve\to0$, from $\cS_0$.
		\item[(M3)] The flow on $\cS_\ve$ converges to the slow flow as $\ve\to0$.
		\item[(M4)] $\cS_\ve$ is $C^k$-smooth.
		\item[(M5)] $\cS_\ve$ is normally hyperbolic and has the same stability properties, i.e., attracting, repelling or saddle-type, with respect to the fast variables as $\cS_0$.
		\item[(M6)] $\cS_\ve$ is usually not unique. In regions that remain at a fixed distance from $\partial\cS_\ve$, all manifolds satisfying items $(M1)-(M5)$ above lie at a Hausdorff distance of order $\cO(\exp(-C/\ve))$ from each other for some $C>0$, $C \in \cO(1)$, as $\ve\to 0$.
	\end{enumerate}
\end{theorem}

In the past, one would refer to Fenichel's theorem as \emph{Geometric Singular Perturbation Theory} (GSPT). As a matter of fact, Fenichel's article~\cite{fenichel1979geometric} is entitled ``Geometric Singular Perturbation Theory for Ordinary Differential Equations''. However, nowadays GSPT includes many more geometric techniques such as the blow-up method and fast-slow normal form theory, among others. A few remarks are now in place.

\newpage

\begin{remark}\leavevmode
	\begin{itemize}[leftmargin=*]
		\item The concept of normal hyperbolicity is applicable, under its appropriate modifications, to invariant manifolds in general and not only to sets of equilibria as in Definition \ref{def:nh}. It turns out, however, that normal hyperbolicity for sets of equilibria is simpler to define. For invariant manifolds in general, normal hyperbolicity is more involved, see \cite{fenichel1971persistence,hirsch1970invariant} and \cite[Chapter 2.3]{kuehn2015multiple}.
		\item Although $\cS_\ve$ is usually not unique, any choice of such manifolds is called \emph{``the''} slow manifold. In turn, one should check if the results one obtains are independent or not on the choice of slow manifold. It is often the case that the results do not depend on such a choice. Furthermore, even though the terms ``critical manifold'' and ``slow manifold'' are sometimes used interchangeably or even equivalently, it is important to make a distinction between them\footnote{In this document $\cC_0$ is a critical manifold and $\cS_0$ is a \emph{compact} subset of $\cC_0$. A perturbation of $\cS_0$, denoted by $\cS_\ve$, is a slow manifold. Note that $\cC_\ve$ may not be well-defined since we do not impose compactness to $\cC_0$; but see~\cite{Eldering1}.}. One of the reasons for this is that, for example, slow manifolds may be extended beyond a non-hyperbolic point, and such extension may not have a relationship with the critical manifold. Many examples of this situation are shown \ch{in Section \ref{sec:summary}}.
		\item In qualitative terms, Fenichel's theorem shows that, under the assumptions of Theorem \ref{thm:fenichel}, a fast-slow system without non-hyperbolic points can be regarded as a regular perturbation of its corresponding singular systems near $\cC_0$.
	\end{itemize}
\end{remark}

The next section describes a geometric technique to locally analyze fast-slow systems around\ch{, a particular class of} non-hyperbolic points.

\subsection{The blow-up method \ch{for fast-slow systems} }\label{sec:blowup}

\ch{In section \ref{sec:bu} we sketched the idea of the blow-up method to desingularize nilpotent singularities of planar vector fields. In this section we provide a review on the blow-up method as is nowadays commonly used for the analysis fast-slow systems with non-hyperbolic \chh{singularities}. To make the relationship between the blow-up method and fast-slow systems clearer, let us first rewrite the $\ve$-family of vector fields \eqref{eq:sfs2} on $\R^{m+n}$ as a single vector field on $\R^{m+n+1}$ of the form
\begin{equation}\label{eq:sfs3}
	\begin{split}
		x' &= f(x,y,\ve)\\
		y' &= \ve g(x,y,\ve)\\
		\ve' &=0.
	\end{split}
\end{equation}
Furthermore, let us assume that the origin $(x,y,\ve)=(0,0,0)$ is an equilibrium point and that $\textnormal D_x f(0,0,0)$ has all its eigenvalues equal to zero (this is indeed the case in the systems discussed in Section \ref{sec:survey}). This means that the origin is a nilpotent singularity of \eqref{eq:sfs3} and, as such, the blow-up method can be adapted to desingularize the origin of \eqref{eq:sfs3}.

\begin{remark}
	It is worth noting that nilpotent singularities are a subset of non-hyperbolic singularities. Thus, not all non-hyperbolic singularities of fast-slow systems may be studied with the blow-up method. However, as we discuss in Section \ref{sec:survey}, there is a large number of fast-slow systems for which, indeed, the origin is nilpotent. 
	In particular, in all fast-slow systems with one-dimensional fast direction ($x\in\R$), a non-hyperbolic singularity is nilpotent. In other cases where the singularity is non-hyperbolic but not nilpotent, a preliminary transformation may bring a fast-slow system into a suitable form to be analyzed via the blow-up method. 
	For example, a centre manifold reduction may lead to a reduced fast-slow system with a nilpotent singularity~\cite{dumortier1991local}. 
	In some other cases, as described in Section \ref{sec:delayedHopf}, a change of coordinates brings the system into an appropriate form.	
\end{remark}

}%
Although there are several (equivalent) versions and improvements of the blow-up method, we restrict to the quasihomogeneous case as it is more commonly used nowadays. For further information see \cite[Chapter 7]{kuehn2015multiple} and references therein.

Let $X:\R^{m+n+1}\to\R^{m+n+1}$ be the vector field, which in coordinates is defined by \eqref{eq:sfs3}, and let $\bbS^N$ denote the $N$-th dimensional sphere. To use spheres, cylinders, or related spaces as blown-up spaces is often very convenient, yet not necessary~\cite{KuehnHyp}. Next, we can formally define the blow-up transformation most commonly used in fast-slow systems:

\begin{definition}[Quasihomogeneous blow-up] Consider a vector field $X:\R^{m+n+1}\to\R^{m+n+1}$ defined by \eqref{eq:sfs3} and assume that $X(0)=0$. Let $\alpha=(\alpha_1,\ldots,\alpha_m)\in\bbN^m_0$, $\beta=(\beta_1,\ldots,\beta_m)\in\bbN^n_0$ and $\gamma\in\bbN_0$. Let the generalized polar transformation $\phi:\bbS^{m+n}\times I\to\R^{m+n+1}$ be defined by
	\begin{equation}\label{eq:phi}
		\phi(\bar x,\bar y,\bar \ve,r) = (r^\alpha\bar x, r^\beta\bar y,r^\gamma\bar\ve)=(x,y,\ve),
	\end{equation}
where $(\bar x, \bar y, \bar\ve)=(\bar x_1,\ldots,\bar x_m,\bar y_1,\ldots,\bar y_n,\bar\ve)\in\bbS^{m+n}$, $r\in I$, and $I\subseteq\R$ is an interval containing the origin. Here we use the multi-index notation $r^\alpha\bx=(r^{\alpha_1}\bar x_1,\ldots,r^{\alpha_m}\bar x_m)$, and similarly for $r^\beta\bar y$. The \emph{quasihomogeneous} blow-up of the vector field $X$, denoted as $\bar X$, is defined by
\begin{equation}\label{eq:bX}
	\bar X = \D\phi^{-1}|_{(\bar x,\bar y,\bar\ve,r)}\circ X \circ \phi(\bar x, \bar y,\bar\ve,r).
\end{equation}
\end{definition}
We note that $\phi$ maps the sphere $\cB_0\coloneqq \mathbb{S}^{m+n}\times\left\{ 0 \right\}$ to the origin in $\R^{m+n+1}$, while $\phi^{-1}$ maps $0\in\R^{m+n+1}$ to $\cB_0$. Hence, the operation $\phi^{-1}$ is called (quasihomogeneous) blow-up while $\phi$ is called (quasihomogeneous) blow-down. The word quasihomogeneous reflects the fact that the exponents appearing in \eqref{eq:phi} are not necessarily the same. We omit the term ``quasihomogeneous'' when all exponents $(\alpha,\beta,\gamma)$ are equal to $1$.

It follows from \eqref{eq:bX} that $\bar X$ and $X$ are conjugate for $r>0$, meaning that there exists a one-to-one mapping between trajectories of $X$ and trajectories of $\bar X$ outside $\cB_0$. Moreover, it can be shown that $\bar X$ is well defined at $r=0$ \cite{kuehn2015multiple}. Due to the presence of non-hyperbolic singularities, and depending on the choice of the exponents, it is usually the case that the system denoted by $\bar X$ vanishes on $\cB_0$. In fact, let $j_\ell(X)$ denote the $\ell$-jet~\cite{golubitsky2012stable} of $X$ at the origin. If $j_\ell(X)=0$ for $\ell=0,1,\ldots,k$ and $j_{k+1}(X)\neq 0$, then we define the \emph{desingularized vector field} $\tilde X = \frac{1}{r^k}\bar X$. Now $\tilde X$ does not vanish at $\cB_0$.  Since $\bar X$ and $\tilde X$ are smoothly equivalent for $r>0$, all the information obtained from $\tilde X$ is equivalent to that of $\bar X$ outside $\cB_0$. However, since $\tilde X$ does not vanish any more along $\left\{ r=0 \right\}$, we may try to infer the dynamics of $\tilde X$ for $r>0$ small from the restriction $\tilde X|_{\left\{ r=0 \right\}}$. This greatly simplifies the analysis, since usually we find that $\tilde X$ has semi-hyperbolic singularities, hyperbolic singularities, or no singularities at all. Finally, due to the
equivalences between $X,\bar X$, and $\tilde X$, we conclude that the flow of
$\tilde X$ for $r>0$ sufficiently small provides a complete description of the flow
of $X$ for $\ve>0$ sufficiently small.

When we study high dimensional problems, say for $m+n>2$, working with polar coordinates can become cumbersome. Then, we rather work in charts that cover the blow-up space. In each of the charts we can define local coordinates and a corresponding \emph{local} vector field. In practice, what we do to define local coordinates in a chart is to fix one of the blow-up coordinates to $\pm 1$. This approach is called ``directional blow-up''. For example, to perform a blow-up in the $\bar\ve$-direction we would define new coordinates according to $\phi:\R^{n+m+1}\to\R^{n+m+1}$ given by
\begin{equation}\label{eq:phi_central}
		\phi(\bar x,\bar y,\bar \ve,r) = (r^\alpha\bar x, r^\beta\bar y,r^\gamma),
\end{equation}
that is by fixing $\bar\ve=1$. Similarly, we can define blow-ups in any of the other directions.
\begin{remark}
The chart $K\coloneqq \left\{ \bar\ve=1 \right\}$ is the most important one and it is called \emph{the rescaling chart}, \emph{the family chart} or \emph{the central chart}. The rest of the charts are often referred to as \emph{phase-directional charts}.
\end{remark}

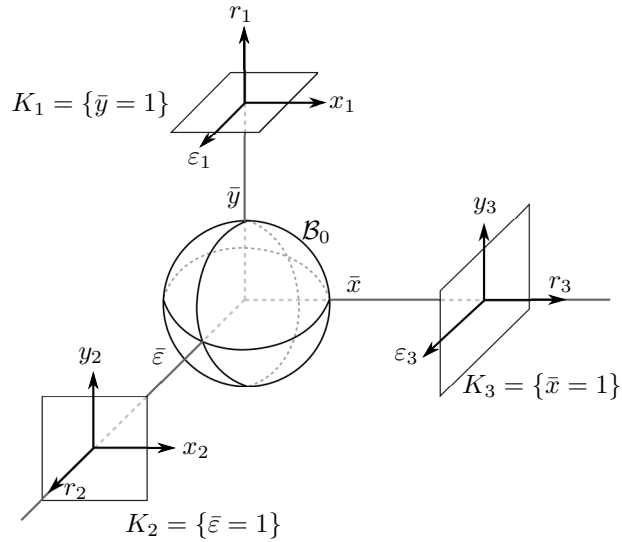
\begin{figure}[htbp]\centering
\begin{tikzpicture}
	\node at (0,0){
	\includegraphics{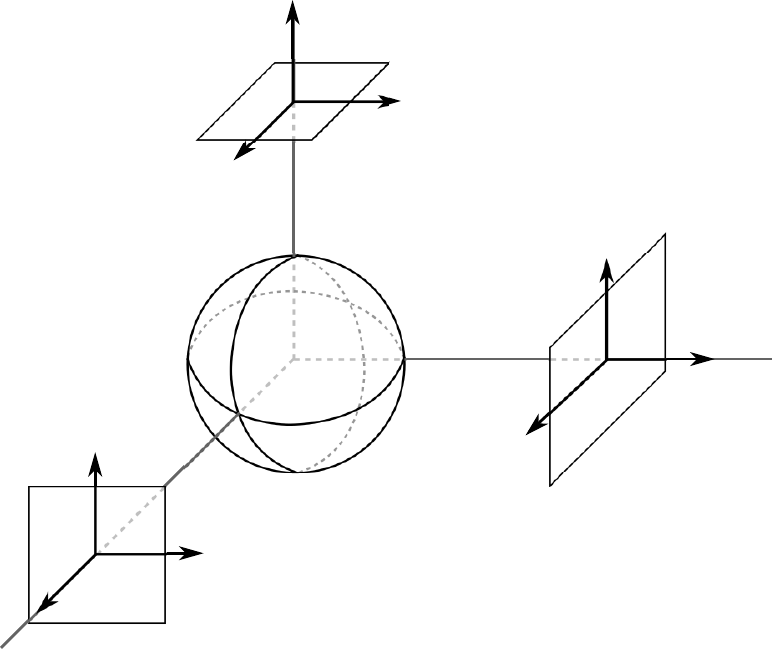}
	};
	\node at (0.5,-0.15) {$\bar x$};
	\node at (-1.1,1) {$\bar y$};
	\node at (-2.1,-1.1) {$\bar \ve$};
	\node at (-1.55,1.5) {$\ve_1$};
	\node at (0.35,2.25) {$x_1$};
	\node at (-1,3.45) {$r_1$};
	\node at (-3.2,-2.9) {$r_2$};
	\node at (-1.6,-2.35) {$x_2$};
	\node at (-3,-1.1) {$y_2$};
	\node at (3.25,-.15) {$r_3$};
	\node at (2.25,.9) {$y_3$};
	\node at (1.2,-1.1) {$\ve_3$};
	\node at (-3.,2.25){$K_1=\left\{ \bar y =1\right\}$};
	\node at (-1.5,-3.35){$K_2=\left\{ \bar \ve =1\right\}$};
	\node at (3.,-1.5){$K_3=\left\{ \bar x =1\right\}$};
	\node at (0,0.55) {$\cB_0$};
\end{tikzpicture}
	\caption{Sketch of the blown-up space and of some of the directional charts. In practice, via the blow-up method, we study local vector fields defined in the charts, and then ``glue'' trajectories and other invariant objects together to describe the dynamics in a small neighborhood of $\cB_0$, which in turn provides the dynamics of a fast-slow system around the origin for $\ve>0$ sufficiently small.}
	\label{fig:blowup}
\end{figure}

Directional blow-ups induce local vector fields on each of the (Euclidean) charts. Once the analysis of the relevant local vector fields is performed, one can overlap suitable regions of the charts and match the flow on such charts via the so-called \emph{matching maps} (or \emph{transition maps}) to describe the dynamics all around $\bbS^{m+n}\times I$. In particular, this process allows us to track invariant objects, principally centre manifolds~\cite{carr2012applications}, across the blow-up space. A schematic representation of the blow-up map is provided in Figure \ref{fig:blowup}.

In summary, through the blow-up method one attempts to make the analysis of fast-slow systems near non-hyperbolic equilibrium points more accessible. Yet, the analysis of the vector fields obtained after blow-up is still highly non-trivial due to several difficult obstacles such as the nonlinearity of the local vector fields, presence of semi-hyperbolic equilibria, and the appearance of resonances, just to mention a few. It may even happen that, in one of the charts, there is a singularity that is still \ch{nilpotent}, in which case the blow-up procedure may be applied again. Furthermore, obtaining the right weights $(\alpha,\beta,\gamma)$ is not always trivial. One way in which one can obtain such weights is via the Newton Polyhedron~\cite{alvarez2011survey,dumortier1977singularities,kuehn2015multiple}.

In the next section we are going to review some of the most fundamental and (mathematically) influential works that have established the blow-up method as one preferred technique to study fast-slow systems in the neighborhood of non-hyperbolic points.

\section{Survey}\label{sec:survey}

Due to their importance to the theory, we first focus on fold points, folded singularities and canards. Later on, we briefly discuss instances where the blow-up method has been used to study singularities beyond the fold such as Hopf, transcritical, pitchfork, cusp, and Bogdanov-Takens singularities.

\subsection{The generic fold}\label{sec:genericfold}

The generic fold is the first natural situation, where a non-hyperbolic point appears in planar fast-slow systems. The exposition of this section mainly summarizes the results of~\cite{krupa2001extending}.

\begin{definition}[Generic fold point] Consider \eqref{eq:sfs1} with $(x,y)\in\R^2$. We say that $(x_0,y_0)$ is a \emph{generic fold point} if
\begin{equation}
	f(x_0,y_0,0)=0, \qquad \frac{\partial f}{\partial x}(x_0,y_0,0)=0,
\end{equation}
hold together with the non-degeneracy conditions
\begin{equation}\label{eq:nodeg}
	\frac{\partial^2 f}{\partial x^2}(x_0,y_0,0)\neq0, \qquad \frac{\partial f}{\partial y}(x_0,y_0,0)\neq0, \qquad g(x_0,y_0,0)\neq0.
\end{equation}
\end{definition}

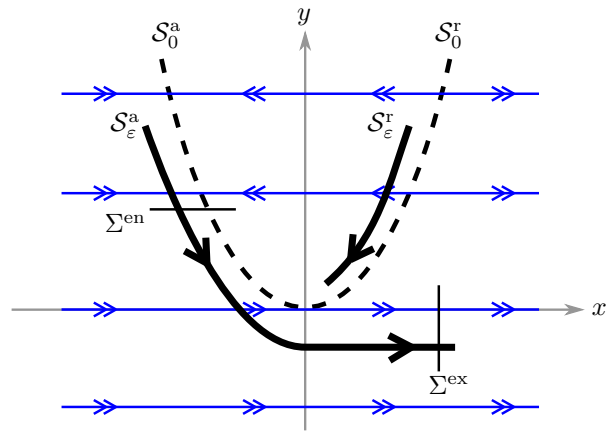
\begin{figure}[htbp]\centering
\begin{tikzpicture}
	\node at (0,0){
	\includegraphics[scale=1.25]{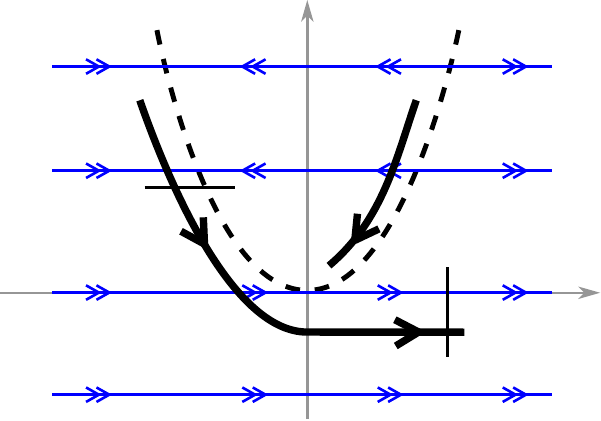}
	};
	\node at (4,-1.05) {$x$};
	\node at (0.075,2.85) {$y$};
	\node at (-1.75,2.6) {$\cS_0^\txta$};
	\node at ( 2,2.6) {$\cS_0^\txtr$};
	\node at (-2.3,1.4) {$\cS_\ve^\txta$};
	\node at ( 1.1,1.4) {$\cS_\ve^\txtr$};
	\node at (-2.3,.1) {$\Sigma^{\EN}$};
	\node at (2,-2.05) {$\Sigma^{\EX}$};
\end{tikzpicture}
\caption{Schematic of a fast-slow system near a generic fold point. Up to leading order terms, the critical manifold $\cC_0$, shown in dashed, is given by $\cC_0=\left\{ (x,y)\in\R^2\, | \, y=x^2 \right\}$. The (blue) lines with double arrows depict the dynamics of the layer equation. Thus, $\cC_0$ has an attracting ($\cS_0^\txta$) and a repelling branch ($\cS_0^\txtr$). Away from the fold point, Fenichel's theorem shows that, for $\ve>0$ sufficiently small, $\cS_0^\txta$ and $\cS_0^\txtr$ are smoothly perturbed to invariant manifolds (in this case trajectories)  $\cS_\ve^\txta$ and $\cS_\ve^\txtr$ respectively. The analysis (via the blow-up method) shows that $\cS_\ve^\txta$ can be extended beyond the fold point as depicted in the figure. In particular, one can show that the distance between the $x$-axis and the intersection  $\cS_\ve^\txta\cap\Sigma^{\EX}$ is of order $\cO(\ve^{2/3})$.  }
\label{fig:fold1}
\end{figure}

Without loss of generality, we assume that $(x_0,y_0)=(0,0)$ and we may choose signs of \eqref{eq:nodeg}. Thus a fast-slow system with a generic fold point at the origin can be given as the (local) canonical form~\cite{krupa2001extending}
\begin{equation}\label{eq:fold}
\begin{split}
		x' &= -y + x^2 + \cO(\ve,xy,y^2,x^3),\\
		y' &= \ve (-1 + \cO(x,y,\ve)).
\end{split}
\end{equation}
We observe that the corresponding critical manifold
\begin{equation}
	\cC_0 = \left\{ (x,y)\in\R^2\, | \, -y + x^2 + \cO(\ve,xy,y^2,x^3)=0 \right\}
\end{equation}%
is locally a parabola as shown in Figure \ref{fig:fold1}.

Let $\cS_0^\txta$ (resp. $\cS_0^\txtr$) denote a compact subset of the attracting (resp.~repelling) section of the critical manifold $\cC_0$. It is known from Fenichel's theorem that a slow manifold $\cS_\ve^\txta$ (resp.~$\cS_\ve^\txtr$) exists as a smooth perturbation of $\cS_0^\txta$ (resp.~$\cS_0^\txtr$). However, Fenichel's theorem only guarantees the existence of such manifolds away from the fold point. The main goal of \cite{krupa2001extending} is to extend $\cS_\ve^\txta$ beyond the fold point.
Let $\rho>0$ be small and $I\subset\R$ be a suitable interval, and define the sections
\begin{equation}
	\begin{split}
		\Sigma^{\EN} &= \left\{ (x,y)\in\R^2\, | \, x\in I, \, y=\rho^2 \right\},\\
		\Sigma^{\EX} &= \left\{ (x,y)\in\R^2\, | \, x=\rho, \, y\in\R \right\}.
	\end{split}
\end{equation}

Here $I$ being a suitable interval means that $\Sigma^{\EN}$ intersects $\cS_0^\txta$ but not $\cS_0^\txtr$, see Figure \ref{fig:fold1}. The following theorem states the properties of the slow manifold $\cS_\ve^\txta$ as it passes through the origin (see Figure \ref{fig:fold1} for a sketch of the result).

\begin{theorem}[{\cite[Theorem 2.1]{krupa2001extending}}]\label{thm:fold} Let $\Pi:\Sigma^{\EN}\to\Sigma^{\EX}$ be the transition map for the flow of \eqref{eq:fold}. Then, there exists $\ve_0>0$ such that the following assertions hold for $\ve\in(0,\ve_0]$:
\begin{enumerate}[leftmargin=*]
	\item[(F1)] The manifold $\cS_\ve^\txta$ passes through $\Sigma^{\EX}$ at a point $(\rho,h(\ve))$, where $h(\ve)\in \cO(\ve^{2/3})$.
	\item[(F2)] The transition $\Pi$ is a contraction with contraction rate $\cO(\exp(-C/\ve))$, where $C>0$.
\end{enumerate}
\end{theorem}

We note the importance of the two \ch{statements (F1) and (F2) contained Theorem \ref{thm:fold} }. First of all, the estimate $h(\ve)\in \cO(\ve^{2/3})$ indicates that $\cS_\ve^\txta$ is not a regular perturbation of the singular set $\cS_0^\txta\cup\left\{ (x,0)\in\R^2 \, | \, x\geq 0 \right\}$. In other words, it is not possible to consider an asymptotic series in positive integer weights of $\ve$ as an expansion of $\cS_\ve^\txta$ in a small neighborhood of the fold point. Next, the second \ch{statement shows} that the exponential contraction towards $\cS_\ve^\txta$, which before the fold point is guaranteed by Fenichel's theorem, is not lost after crossing the fold point. We observe that this fact is not a regular perturbation argument derived from the corresponding layer equation, but can be proven via the blow-up analysis.

Let us now point-out the key steps of the proof of Theorem \ref{thm:fold}. The blow-up is defined by
\begin{equation}
\label{eq:foldbu1}
	x=\br\bx,\qquad y=\br^2\by,\qquad\ve=\br^3\be.
\end{equation}
The weights of the blow-up map follow formally from the quasi-homogeneity~\cite{dumortier1991local} of the function $-y+x^2$. On a practical level, one can also guess good weights by observing that the scaling $(x,y,t)\mapsto (\ve^{1/3}x,\ve^{2/3}y,\ve^{-1/3}t)$ removes $\ve$ to leading-order from~\eqref{eq:fold}. Using~\eqref{eq:foldbu1}, the blow-up analysis is carried over three charts, namely
\begin{equation}
	K_1=\left\{ \by=1 \right\}, \quad K_2=\left\{ \be=1 \right\}, \quad K_3=\left\{ \bx=1 \right\}.
\end{equation}
The behaviour of the trajectories passing through $\Sigma^{\EN}$ as they approach (but do not cross) a small neighbourhood of the origin is studied in chart $K_1$. The main objects of study in this chart are $2$-dimensional centre manifolds associated to semi-hyperbolic equilibrium points. To be more precise, in chart $K_1$ one can show the existence of a $2$-dimensional attracting and a $2$-dimensional repelling centre manifolds that are denoted by $\cM_1^\txta$ and $\cM_1^\txtr$ respectively. In fact, sections of these centre manifolds correspond to perturbations of $\cS_0^\txta$ and $\cS_0^\txtr$, that is to slow manifolds, along normally hyperbolic parts of the critical manifold. The advantage of using the blow-up method is that it allows to track invariant objects (in this case centre manifolds) as they pass through all the charts.

Next, in chart $K_2$, we study the behaviour of trajectories in a small neighbourhood of the origin. The main object of study in this chart is a Riccati differential equation~\cite[pp.~68-72]{mishchenko1980differential} of the form
\begin{equation}\label{eq:riccati}
	\begin{split}
		x_2' &= -y_2 + x_2^2,\\
		y_2' &= -1,\\
		r_2' &= 0
	\end{split}
\end{equation}
where we use the subscript $2$ to indicate local coordinates in the chart $K_2$. A detailed analysis of \eqref{eq:riccati} allows, in particular, to extend $\cM_1^\txta$ found in chart $K_1$ to an invariant manifold $\cM_2^\txta$ beyond the fold point. The previous is due to the fact that sections of the centre manifold $\cM_1^\txta$ (in chart $K_1$) can be matched with solutions of \eqref{eq:riccati} with $r_2$ constant (in chart $K_2$).
Afterwards, in chart $K_3$, we study the dynamics of the trajectories, as they leave a small neighbourhood of the fold point. More specifically, similar to the situation in the previous two charts, via a change of coordinates we can define an invariant object $\cM_3^\txta$ corresponding to $\cM_2^\txta$ but defined in the local coordinates of chart $K_3$. Then we can pay special attention to trajectories on $\cM_3^\txta$. The main object of study in this chart is a $3$-dimensional nonlinear ODE with a resonant hyperbolic saddle at the origin\footnote{We recall the definition of resonant equilibrium point: consider an ODE of the form $\dot x=Ax$, where $x\in\R^n$, and let $A$ have all its eigenvalues $\lambda_1,\ldots,\lambda_n$ with non-zero real part. We say that $x=0$ is resonant if there is a relation $\lambda_j=\sum_{i=1}^n m_i\lambda_i$, where the $m_i$'s are non-negative integers with $\sum_{i=1}^n m_i\geq 2$.}. The fact that the origin is resonant is a major obstacle as it prevents to linearize the nonlinear system in order to have a detailed description of the flow near the origin. In fact, dealing with the aforementioned resonance provides that the next term in the expansion of $h$ is of the order $\cO(\ve\ln\ve)$. Once the analysis in these three charts is performed, the transition $\Pi$ can be described by overlapping the individual transitions in each of the charts via the corresponding matching maps.

We emphasize that the order of the function $h$ was already known in \cite{mishchenko1980differential}. In fact, rigorous asymptotic expansions of the function $h$ appear in \cite[equation (16.11) page 106]{mishchenko1980differential}. It is important to note that, although during the 80's the blow-up method was not known per s\'e in the context of fast-slow systems, the analysis performed in \cite{mishchenko1980differential} starts from a rescaling similar to the blow-up presented above, see (8.1) and (8.2) in \cite{mishchenko1980differential}. More recently, in \cite[Theorem 2]{van2005asymptotic} a rigorous asymptotic expansion of the function $h$ is obtained by combining the blow-up method with matched asymptotic analysis. The main idea of \cite{van2005asymptotic} is to compute, in each chart, asymptotic expansions of the centre manifolds $\cM_i^\txta$ mentioned above, and then match such local expansions across the blow-up space.

\subsection{Planar canards}\label{sec:planarcanards}

\begin{figure}[htbp]\centering
	\begin{tikzpicture}
	\node at (0,0) {\includegraphics[scale=1.25]{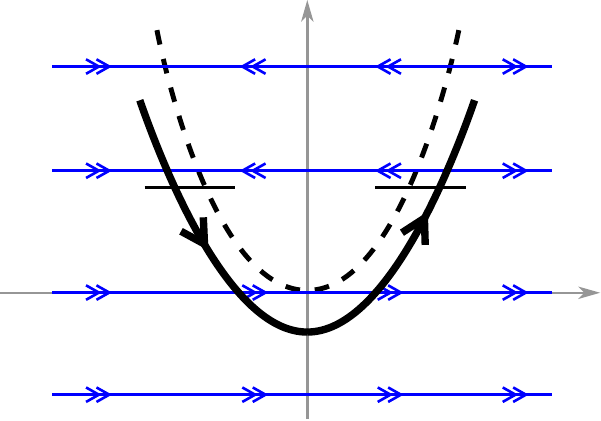}};
	\node at (4,-1.05) {$x$};
	\node at (0.075,2.85) {$y$};
	\node at (-1.75,2.6) {$\cS_0^\txta$};
	\node at ( 2,2.6) {$\cS_0^\txtr$};
	\node at (-2.3,1.4) {$\cS_\ve^\txta$};
	\node at ( 2.5,1.4) {$\cS_\ve^\txtr$};
	\node at (-2.3,.1) {$\Sigma^{\EN}$};
	\node at ( 2.4,.1) {$\Sigma^{\EX}$};
	\end{tikzpicture}
	\caption{Schematic representation of a maximal canard. Note that $\cS_\ve^\txta$ extends precisely to $\cS_\ve^\txtr$ beyond the origin. Trajectories passing  through $\Sigma^\EN$ exponentially close to $\cS_\ve^\txta$ follow  $\cS_\ve^\txtr$ for time of order $\cO(1)$ before being repelled from it.}
	\label{fig:canard}
\end{figure}
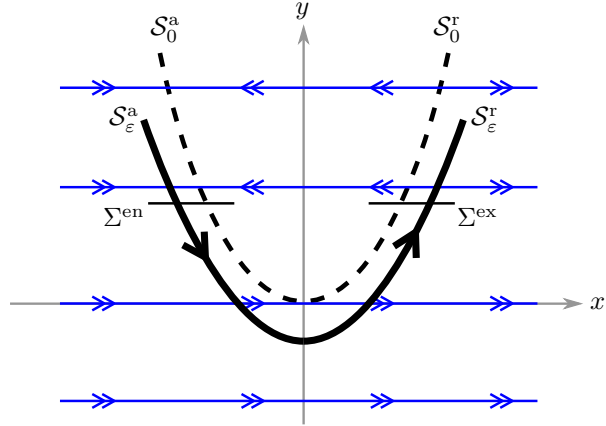

In Section \ref{sec:genericfold} we have seen that trajectories near a generic fold point first follow the attracting part of the critical manifold, and then follow the fast direction. Perhaps counter-intuitively, for certain one-parameter families of planar fast-slow systems, there are trajectories  near a fold singularity that can closely follow the unstable region of the critical manifold for time of order $\cO(1)$. Such type of trajectories are called \emph{canards}, see Figure \ref{fig:canard} (and the rest of the figures in this Section). In this section we summarize the results of \cite{dumortier1996canard,krupa2001relaxation} and Section 3 of \cite{krupa2001extending} dealing with the description of canards in planar fast-slow systems.

\begin{remark}
The blow-up technique was first introduced for the analysis of fast-slow systems in \cite{dumortier1996canard} to explain, in a geometric way, the canard phenomenon in the van der Pol oscillator. Canards were already explained in \cite{benoit1981chasse} with the language of non-standard analysis. Although \cite{dumortier1996canard} was the first work describing geometrically the canard phenomenon, we start from \cite{krupa2001extending} as it follows several arguments analogous to the generic fold point.
\end{remark}

We start by considering a one-parameter family of planar fast-slow systems given by
\begin{equation}\label{eq:sfs-canard}
	\begin{split}
		x' &= f(x,y,\omega,\ve),\\
		y' &= \ve g(x,y,\omega,\ve),
	\end{split}
\end{equation}
where $\omega\in\R$. Just as in the generic fold case of Section \ref{sec:genericfold} we may assume that a fold point exists at the origin, say for $\omega=0$.

\begin{definition}[Canard point] Consider \eqref{eq:sfs-canard}. We say that $(x,y)=(0,0)$ is a \emph{canard point} if
\begin{equation}\label{eq:cp1}
 	f(0,0,0,0)=0, \quad \frac{\partial f}{\partial x}(0,0,0,0)=0, \quad g(0,0,0,0)=0
\end{equation}
	hold together with the non-degeneracy conditions
\begin{equation}\label{eq:cp2}
	\frac{\partial^2f}{\partial x^2}(0,0,0,0)\neq0, \quad \frac{\partial f}{\partial y}(0,0,0,0)\neq0, \quad \frac{\partial g}{\partial x}(0,0,0,0)\neq0, \quad \frac{\partial g}{\partial \omega}(0,0,0,0)\neq0.
\end{equation}
\end{definition}
The main difference with the generic fold case of Section \ref{sec:genericfold} is that there is an equilibrium point at the origin, which is imposed by $g(0,0,0,0)=0$. The non-degeneracy conditions imply a transverse intersection at the origin of the nullclines of \eqref{eq:sfs-canard}. The critical manifold of \eqref{eq:sfs-canard} is, up to leading order terms, the same as in Section \ref{sec:genericfold}. So, we use the same notation and do not repeat its definition here.
Similarly to the generic fold case, it can be shown that a planar fast-slow system has a local expression near a canard point given by
\begin{equation}\label{eq:canard-nf}
	\begin{split}
		x' &= -y h_1(x,y,\omega,\ve) + x^2 h_2(x,y,\omega,\ve) + \ve h_3(x,y,\omega,\ve),\\
		y' &= \ve\left( x h_4(x,y,\omega,\ve) - \omega h_5(x,y,\omega,\ve) + yh_6(x,y,\omega,\ve) \right),
	\end{split}
\end{equation}
where
\begin{equation}
	\begin{split}
		h_3 &= \cO(x,y,\omega,\ve),\\
		h_j &= 1+\cO(x,y,\omega,\ve), \quad j=1,2,4,5.
	\end{split}
\end{equation}

\begin{remark}
	Up to leading order terms \eqref{eq:canard-nf} reads as
	\begin{equation}
		\begin{split}
			x' &= -y+x^2,\\
			y' &= \ve(x-\omega).
		\end{split}
	\end{equation}
	Let, for a moment, $\ve > 0$. Then a Hopf bifurcation occurs at $\omega=0$. Moreover, the norm of the corresponding eigenvalues tends to zero (and to infinity in the slow time scale) as $\ve\to0$. This situation is known as \emph{singular Hopf bifurcation}~\cite{braaksma1998singular,baer1986singular}.
\end{remark}

Define the constants
\begin{equation}
\begin{split}
	a_1 &= \frac{\partial h_3}{\partial x}(0,0,0,0),\\
	a_5 &= h_6(0,0,0,0), \\
	A &= -\frac{\partial h_1}{\partial x}(0,0,0,0)+3\frac{\partial h_2}{\partial x}(0,0,0,0)-2\frac{\partial h_4}{\partial x}(0,0,0,0)-2a_5
\end{split}
\end{equation}
and the sections
\begin{equation}
	\Sigma^j = \left\{ (x,y)\in\R^2\,|\, x\in I_j,\, y=\rho^2 \right\}, \qquad \qquad j=\EN,\EX,
\end{equation}
where $\rho>0$ and $I_j$ is a suitable small interval, see Figure \ref{fig:canard}. Let $q_{\EN,\ve}=\Sigma^{\EN}\cap\cS_{\ve}^\txta$, $q_{\EX,\ve}=\Sigma^{\EX}\cap\cS_{\ve}^\txtr$ and denote by $\Pi:\Sigma^\EN\to\Sigma^\EX$ the transition induced by \eqref{eq:sfs-canard}. The main result that describes the behavior of the slow manifolds $\cS_{\ve}^\txta$ and $\cS_{\ve}^\txtr$ around the canard point is as follows.

\begin{theorem}[{\cite[Theorem 3.1]{krupa2001extending}}]\label{thm:canard1} Assume that \eqref{eq:sfs-canard} has a canard point at the origin, i.e. satisfies \eqref{eq:cp1} and \eqref{eq:cp2}. Assume that the solution $x_0(t)$ of the reduced problem connects $\cS_0^\txta$ to $\cS_0^\txtr$. Then there exists $\ve_0>0$ and a smooth function $\omega_c(\sqrt\ve)$ defined on $[0,\ve_0]$ such that for $\ve\in(0,\ve_0]$ the following assertions hold.
\begin{enumerate}[leftmargin=*]
	\item[(C1)] $\Pi(q_{\EN,\ve})=q_{\EX,\ve}$ if and only if $\omega=\omega_c(\sqrt\ve)$.
	\item[(C2)] The function $\omega_c$ has the expansion
	\begin{equation}
		\omega_c(\sqrt\ve)=-\left( \frac{a_1+a_5}{2} + \frac{1}{8}A \right)\ve+\cO(\ve^{3/2})
	\end{equation}
	\item[(C3)] The transition map $\Pi$ is defined only for $\omega$ in an interval around $\omega_c(\sqrt\ve)$ of size  $\cO(\exp(-C/\ve))$ for some $C>0$.
	\item[(C4)] There is a splitting condition at $\omega=\omega_c(\sqrt\ve)$ given by
	\begin{equation}
		\frac{\partial}{\partial\omega}(\Pi(q_{\EN,\ve})-q_{\EX,\ve})|_{\omega=\omega_c(\sqrt\ve)}>0.
	\end{equation}
\end{enumerate}
\end{theorem}

A schematic representation of Theorem \ref{thm:canard1} is given in Figure \ref{fig:canard}. The most important aspect of Theorem \ref{thm:canard1} is that it tells us that for $\ve>0$ there is exactly one value of the parameter $\omega$ for which $\cS_\ve^\txta$ is extended to \emph{exactly} $\cS_\ve^\txtr$, that is, they coincide as stated in (C1). The previous trajectory is usually called \emph{maximal canard} to distinguish it from others that may follow $\cS_\ve^\txtr$ for time of order $\cO(1)$. Another important point of Theorem \ref{thm:canard1} is its third item (C3), which tells us that canard orbits near planar folds occur only for values of $\omega$ within an exponentially small interval. Related to the previous observation we further note that the intersection of $\cS_\ve^\txta$ and $\cS_\ve^\txtr$ is not transversal and thus small $C^1$-perturbations of the parameter $\omega$ destroy the maximal canard; see (C4) and see also Section~\ref{sec:foldedsing}.

Let us now point-out the key ingredients of the analysis to prove Theorem \ref{thm:canard1}. In many ways, the blow-up analysis of a canard point is similar to that of a generic fold point. One of the differences is that now the parameter $\omega$ is included in the blow-up. The appropriate blow-up transformation is defined by
\begin{equation}
	x = \br\bx, \quad y = \br^2\by, \quad \ve=\br^2\be, \quad \omega=\br\bar\omega.
\end{equation}

Due to the definition of the $\Sigma^\EN$ and $\Sigma^\EX$ sections, the blow-up analysis is performed only in two charts, namely $K_1=\left\{ \by=1 \right\}$ and $K_2=\left\{ \be=1 \right\}$. The analysis in chart $K_1$ for the canard point is analogous to the analysis in chart $K_1$ for the generic fold point. The only difference is that in the canard case, the consideration of $\omega$ in the blow-up map increases by one the dimension of the centre manifolds. Thus, in $K_1$, the main object of study are $3$-dimensional centre manifolds $\cM_1^\txta$ and $\cM_1^\txtr$ that are of centre-stable and centre-unstable type respectively. These centre manifolds correspond to invariant manifolds $\cM_2^\txta$ and $\cM_2^\txtr$ in the chart $K_2$. In chart $K_2$ the key object of study is a planar ODE of the form
\begin{equation}\label{eq:K2-canard}
	\begin{split}
		x_2' &= -y_2+x_2^2,\\
		y_2' &= x_2,
	\end{split}
\end{equation}
which is obtained after restriction to $\left\{r_2=\omega_2=0\right\}$; compare with \eqref{eq:riccati}. One may check that \eqref{eq:K2-canard} is completely integrable, having a constant of motion denoted by $H(x_2,y_2)$. This means that the trajectories of \eqref{eq:K2-canard} are determined by the level curves of $H$. Here $(x_2,y_2,r_2,\omega_2)$, as usual, denote local coordinates in $K_2$.

There is a special curve $\gamma_c$ given by the solution of $H=0$. The major importance of $\gamma_c$ is that it precisely connects $\cS_0^\txta$ with $\cS_0^\txtr$. Naturally in the blow up space $\cS_0^\txta$ and $\cS_0^\txtr$ correspond to $\cM_2^\txta|_{r_2=0}$ and $\cM_2^\txtr|_{r_2=0}$ respectively.

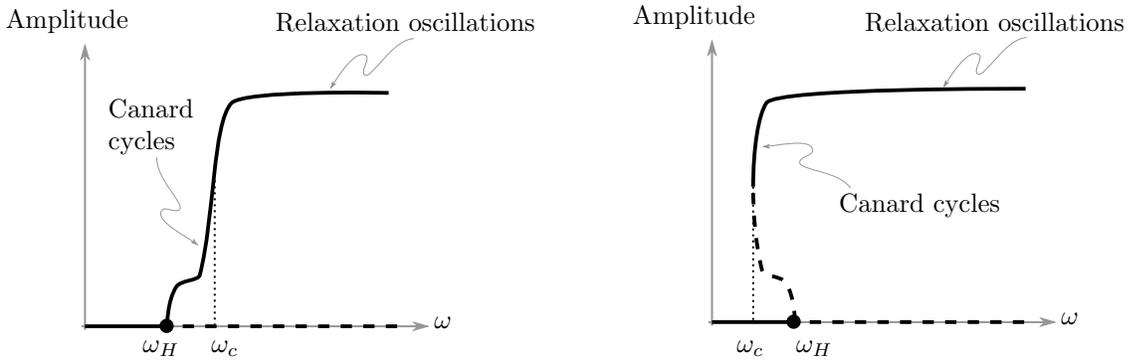
\begin{figure}[htbp]\centering
	\begin{tikzpicture}
		\node at (0,0) {
		\includegraphics{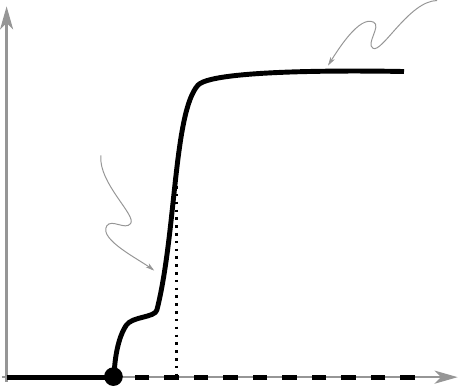}
		};
		\node at (2.5,-1.8) {$\omega$};
		\node at (-2.5,2.2) {Amplitude};
		\node at (-1.25,-2.15) {$\omega_H$};
		\node at (-.4,-2.15) {$\omega_c$};
		\node at (2,2.2) {Relaxation oscillations};
		\node[text width=1.5cm] at (-1.2,.8) {Canard cycles};
	\end{tikzpicture}\hfill
	\begin{tikzpicture}
		\node at (0,0) {
		\includegraphics{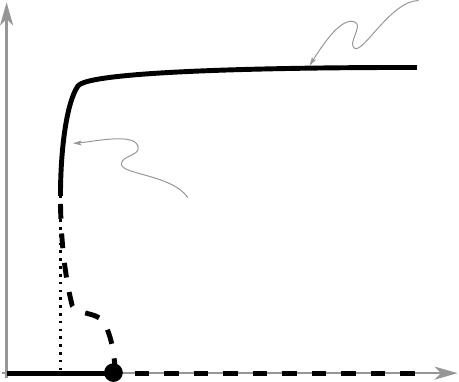}
		};
		\node at (2.5,-1.8) {$\omega$};
		\node at (-2.5,2.2) {Amplitude};
		\node at (-.9,-2.2) {$\omega_H$};
		\node at (-1.75,-2.2) {$\omega_c$};
		\node at (1.5,2.15) {Relaxation oscillations};
		\node at (.5,-.3) {Canard cycles};
	\end{tikzpicture}
	\caption{Bifurcation diagrams corresponding to a canard explosion for $A<0$ on the left and $A>0$ on the right (refer to Theorem \ref{thm:canard1}). A Hopf bifurcation occurs at $\omega=\omega_H$ and a maximal canard when $\omega=\omega_c$. These diagrams depict the amplitude and stability of the limit cycles $\Gamma_{\ve,\omega}$ for fixed $\ve>0$ sufficiently small. In particular, we have that for $A<0$ all canard cycles are stable. }
	\label{fig:canardbifurcation}
\end{figure}

The next step is to investigate how such a connection breaks for $r_2>0,\,\omega_2\neq0$, which will describe how the manifolds $\cM_2^\txta$ and $\cM_2^\txtr$ are connected. For such a purpose an analysis based on the Melknikov method~\cite{guckenheimer2013nonlinear} is used. In essence, the Melknikov method provides an expression for the distance between $\cM_2^\txta$ and $\cM_2^\txtr$, denoted by $d=d(r_2,\omega_2)$. Therefore, for $\ve>0$ the manifolds $\cS_\ve^\txta$ and $\cS_\ve^\txtr$ are connected if and only if $d=0$. Careful estimates of the function $d$ and a blow-down provide the expression for $\omega_c$ in Theorem \ref{thm:canard1} and the breaking argument proving (C4). We finalize the summary of the canard point by mentioning that the role of the constant $A$ has to do with the non-degeneracy of the Hopf bifurcation that occurs for $\omega=0$, see for example \cite[Section 8.3]{kuehn2015multiple}, and Figure \ref{fig:canardbifurcation}.

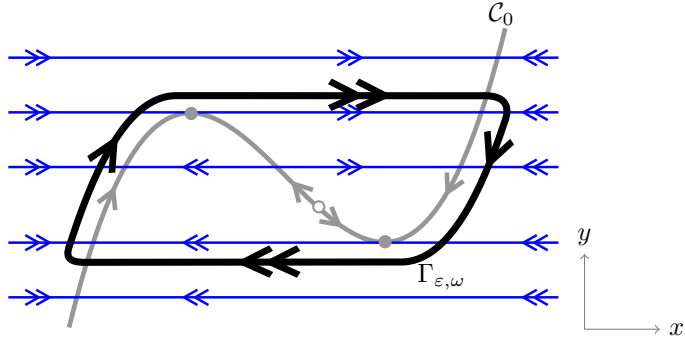
\begin{figure}[htbp]\centering
\begin{tikzpicture}
	\node at (0,0){
		\includegraphics[scale=1.25]{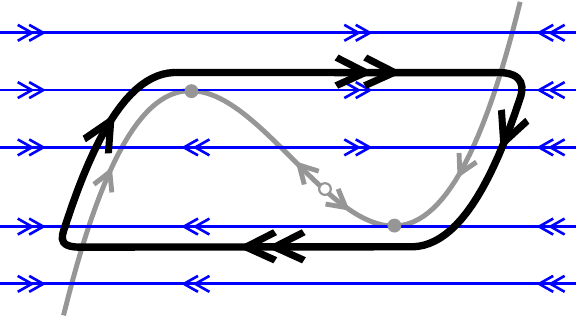}
	};
	\node at (2.9,2.2) {$\cC_0$};
	\node at (2.1,-1.3) {$\Gamma_{\ve,\omega}$};
	\draw[gray,->] (4,-2)--++(1,0) node[right,black]{$x$};
	\draw[gray,->] (4,-2)--++(0,1) node[above,black]{$y$};
\end{tikzpicture}
	\caption{The van der Pol oscillator: the critical manifold is the cubic curve $\cC_0$; the curve in this case is also called S-shaped curve, which is motivated by viewing it in the $(y,x)$-plane. The critical manifold has two fold points indicated by filled circles. The singular limit is represented by the blue lines (layer equation) and by the arrows along $\cC_0$ (CDE). The CDE has an equilibrium point at $x=\omega$ indicated by an empty circle. As it is described in this Section, the parameter $\omega$ is essential to distinguish several types of limit cycles. The curve $\Gamma_{\ve,\omega}$ indicates a typical limit cycle that exists for $\ve>0$ sufficiently small.}
	\label{fig:vdp-sing}
\end{figure}

Next, let us discuss canards in the context of the van der Pol oscillator as they have a particular motivating relevance for the entire subject. Here the main reference is \cite{dumortier1996canard}. To be more precise, \cite{dumortier1996canard} is concerned with \emph{canard cycles}. Canard cycles are closed orbits  having ``a canard piece''. The model of the van der Pol oscillator to be considered is
\begin{equation}\label{eq:vdp}
	\begin{split}
		x' &= y-\frac{x^2}{2}-\frac{x^3}{3},\\
		y' &= \ve(\omega-x).
	\end{split}
\end{equation}
In the singular limit $\ve=0$, the phase portrait is as shown in Figure \ref{fig:vdp-sing} and the critical manifold is a cubic curve given by
\begin{equation}
	\cC_0 = \left\{(x,y)\in\R^2\, |\, y=\frac{x^2}{2}+\frac{x^3}{3} \right\}.
\end{equation}

\begin{figure}[htbp]\centering
	\begin{tikzpicture}
		\node at (0,0) {\includegraphics{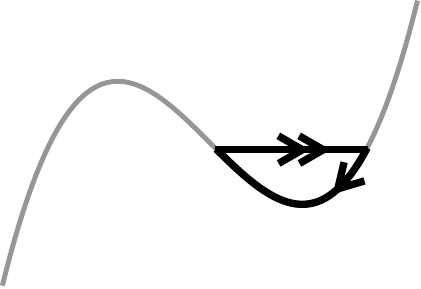}};
		\node at (-5,0) {\includegraphics{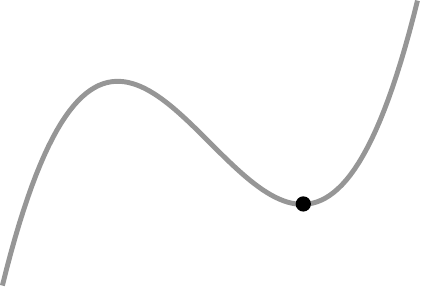}};
		\node at (5,0) {\includegraphics{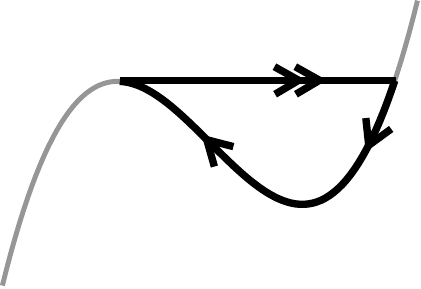}};
		\node at (-2.5,-4) {\includegraphics{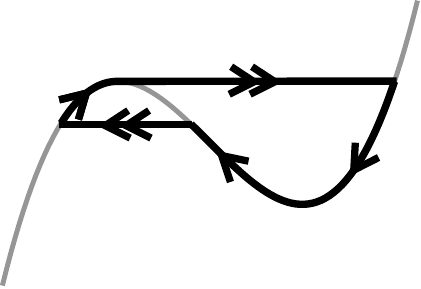}};
		\node at (2.5,-4) {\includegraphics{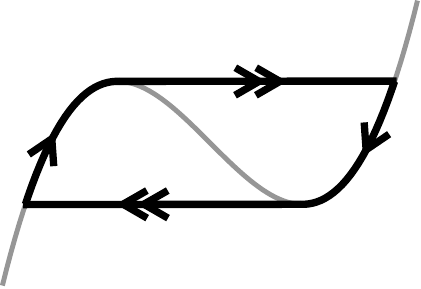}};

		\draw[gray,->] (6,-5.5)--++(1,0) node[right,black]{$x$};
		\draw[gray,->] (6,-5.5)--++(0,1) node[above,black]{$y$};
	\end{tikzpicture}

	\caption{Singular canard cycles (or limit periodic sets) giving rise to canard cycles in the van der Pol oscillator, compare with Figure \ref{fig:vdp-sing}. For $\ve>0$ sufficiently small closed orbits exist as perturbations of these singular cycles.
	}
	\label{fig:canardcycles}
\end{figure}

We note that $\cC_0$ has two fold points namely $(x,y)=\left\{ \left( -1,\frac{1}{6}\right),\,(0,0) \right\}$. Furthermore, there is a canard point at the origin for $\omega=0$. It is well known~\cite{lins1977lienard} that \eqref{eq:vdp} has at most one limit cycle and that when such a limit cycle exists, it is hyperbolic and attracting. Assuming that a limit cycle $\Gamma_{\ve,\omega}$ of \eqref{eq:vdp} exists, the question posed by \cite{dumortier1996canard} can be formulated as follows:  ``what happens to $\Gamma_{\ve,\omega}$ as $\ve\to0$?'' The main contribution of \cite{dumortier1996canard} is the precise asymptotic description of the canard cycles that occur for \eqref{eq:vdp}. In other words, \cite[Theorem 1]{dumortier1996canard} provides asymptotic estimates of the parameter $\omega$ for each of the canard cycles shown in Figure \ref{fig:canardcycles}. All such cycles occur within an exponentially small interval of the parameter $\omega$. This means that the transition between small amplitude canard cycles towards relaxation oscillations occurs for an exponentially small variation of the parameter. This effect is called \emph{canard explosion}.  The analysis presented in \cite{dumortier1996canard} is highly non-trivial: besides introducing the blow-up method for fast-slow systems, it uses techniques of nonlinear analysis such as normal form theory, centre manifold theory, asymptotic analysis, Abelian integrals~\cite{bliss1933algebraic,dumortier2006abelian}, foliations~\cite{camacho2013geometric}, among others.

\begin{figure}[htbp]\centering
	\begin{tikzpicture}
	\node at (0,0){
		\includegraphics[scale=1.25]{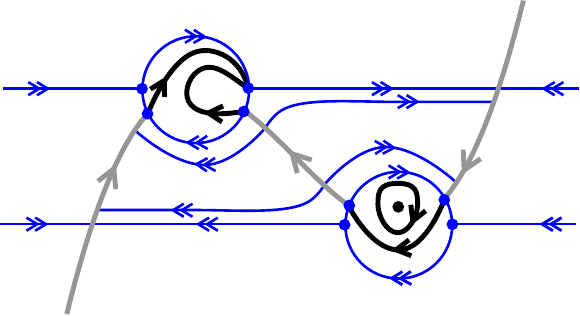}
	};
	\node at (2.6,2.1) {$\cC_0$};
	\draw[gray,->] (5,-2)--++(1,0) node[right,black]{$x$};
	\draw[gray,->] (5,-2)--++(0,1) node[above,black]{$y$};
\end{tikzpicture}
	\caption{Schematic of the families of periodic orbits, parametrized by $\omega$ in \eqref{eq:vdp}, constructed by gluing together trajectories in the blow-up space with trajectories in the singular limit. The flow near the fold and canard points (compare with Figure \ref{fig:vdp-sing}) have been replaced by the flow in the blow-up space. Appropriate combinations of orbits correspond to those depicted in Figure \ref{fig:vdp-sing}. Perturbations of these orbits lead to canard cycles for $\ve>0$ sufficiently small. }
	\label{fig:canardfamilies}
\end{figure}
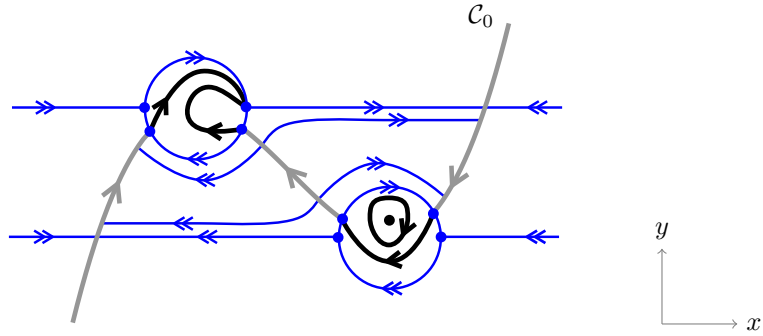

Extensions and generalizations of the results of \cite{dumortier1996canard} can be found in \cite{krupa2001relaxation}. The setting of \cite{krupa2001relaxation} is not confined to the van der Pol oscillator, but it applies to general one parameter families of planar fast-slow systems \eqref{eq:sfs-canard} having and $S$-shaped critical manifold (a particular case being the van der Pol oscillator). Similar to \cite{dumortier1996canard} the main result of \cite{krupa2001relaxation} is the proof that the transition between small amplitude periodic orbits to relaxation oscillations occurs within an exponentially small interval of the parameter $\omega$, the precise statement is in~\cite[Theorem 3.3]{krupa2001relaxation}, while results on the stability of the canard cycles are given in~\cite[Theorems 3.4, 3.5 and 3.6]{krupa2001relaxation}. We briefly summarize the strategy to prove such Theorems as follows. Given an $S$-shaped critical manifold, one assumes that one of the fold points is generic while the other is a canard point. Then, one can locally study the flow near the aforementioned fold points with the techniques reviewed in Section \ref{sec:genericfold} and at the beginning of this Section. Afterwards, families of periodic orbits can be constructed by ``gluing'' together pieces of trajectories in the blow-up space with trajectories of the layer equation and trajectories near normally hyperbolic parts of the critical manifold. A schematic of this construction is shown in Figure \ref{fig:canardfamilies}, while a bifurcation diagram, also representing the canard explosion is shown in Figure \ref{fig:canardbifurcation}.

So far, we have seen that a mechanism leading to canards in planar fast-slow systems is a singular Hopf bifurcation. Some higher codimension bifurcations leading to canards have also been studied via the blow-up method. For example, for the Bogdanov-Takens bifurcation, which is of codimension $2$, we have e.g.~\cite{de2011slow,chiba2011periodic}. The degenerate case, that is a planar system for which normal hyperbolicity is lost up to arbitrary order, has been dealt with in e.g.~\cite{de2006canard}. Another recent discovery is that canards can also be found in cases where the critical manifold is of codimension zero~\cite{kuehn2018duck}\ch{, that is, in the case when the critical manifold is of the same dimension as the phase-space}. Finally, we note that there is strong active research towards proving Hilbert's 16th problem~\cite{smale1998mathematical} based on the tools (and their further generalizations) described in this Section. A non-exhaustive list of references dealing with the cyclicity of canard cycles is: for (fast-slow) Li\'enard systems~\cite{de2008canard,de2011classical,de2011cyclicity,dumortier2007bifurcation,dumortier2009birth,huzak2013limit,huzak2014primary,huzak2016cyclicity,dumortier2007more} and non-Li\'enard systems e.g.~\cite{huzak2018canard}. A proof of finite cyclicity for fast-slow Darboux systems is available in \cite{bobienski2016finite}.

\subsection{Folded singularities and canards in higher dimensions}\label{sec:foldedsing}

In Section \ref{sec:planarcanards} we have seen that canards appear in one-parameter families of planar fast-slow systems having a folded critical manifold. However, canards are degenerate in such a case since they exist only for an exponentially small interval of parameter values and disappear after small perturbations of the parameter. This situation is ``fixed'' in fast slow-systems with at least two slow variables, where canards are generic.

Canards in $\R^3$ were first studied in \cite{benoit1983systemes,benoit1990canards,benoit2001perturbation} with techniques from non-standard analysis; see also \cite{bobkova2002duck} for a different approach in which the critical manifold is assumed to have a generic self-intersection. In this section we review some of the works that, based on the blow-up method, have provided a geometric description of canards in fast-slow systems with two, or more, slow variables.

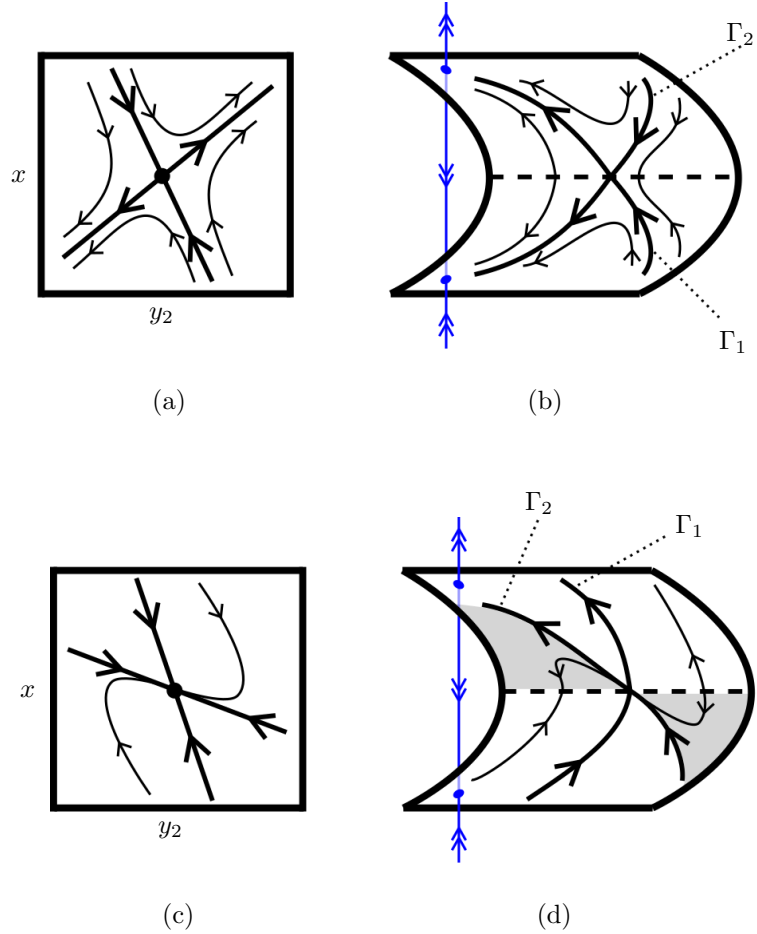
\begin{figure}[htbp]\centering
\begin{tikzpicture}
	\node at (0,0){
	\includegraphics[scale=1.25]{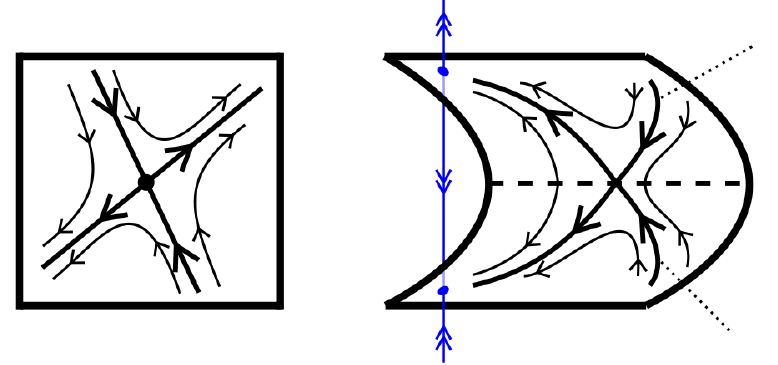}
	};
	\node at (4.5,-2.2) {$\Gamma_1$};
	\node at (4.6, 1.9) {$\Gamma_2$};
	\node at (-3.1,-1.85) {$y_2$};
	\node at (-5, 0) {$x$};
	\node at (-3,-3) {(a)};
	\node at ( 2,-3) {(b)};
	\end{tikzpicture}

	\vspace*{.75cm}

	\begin{tikzpicture}
	\node at (0,0){
	\includegraphics[scale=1.25]{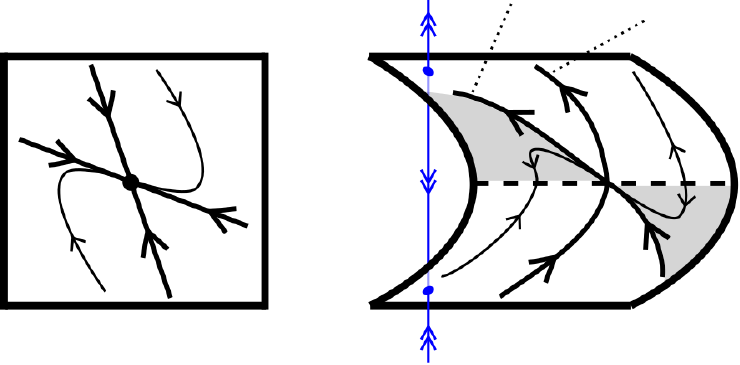}
	};
	\node at (3.8,2.2) {$\Gamma_1$};
	\node at (1.8,2.5) {$\Gamma_2$};
	\node at (-3.1,-1.85) {$y_2$};
	\node at (-5, 0) {$x$};
	\node at (-3,-3) {(c)};
	\node at ( 2,-3) {(d)};
	\end{tikzpicture}
	\caption{ Phase-portraits for the folded saddle (first row) and folded node (second row) cases of Table \ref{table:class}. The planar phase-portraits (a),(c) correspond to the flow of the desingularized vector field \eqref{eq:desing}, while the phase-portraits (b),(d) are the corresponding reduced flows on $\cC_0$. Note that trajectories on the stable (i.e.~lower) part of $\cC_0$ are given by trajectories of the desingularized system for $x<0$ while the trajectories on the unstable (i.e.~upper) part of $\cC_0$ are given by reversing the direction of the trajectories of the desingularized system for $x>0$. In the saddle case (a),(b) the trajectory $\Gamma_1$ is a true canard while $\Gamma_2$ is a faux canard. Note that these are the only singular canard trajectories. In the node case (c),(d) $\Gamma_1$ is a singular canard  tangent to the strong eigendirection while $\Gamma_2$ is a singular canard  tangent to the weak eigendirection. Observe that, in contrast to the saddle case, in the node case we have many other singular canard trajectories indicated by the shaded region in $\cC_0$. }
	\label{fig:desing-example}
\end{figure}

Let us start with a fast-slow system with two slow and one fast variable given by
\begin{equation}\label{eq:sfs3d}
	\begin{split}
		x'  &= f(x,y_1,y_2,\ve),\\
		y_1' &= \ve g_1(x,y_1,y_2,\ve),\\
		y_2' &= \ve g_2(x,y_1,y_2,\ve),
	\end{split}
\end{equation}
where $(x,y)=(x,y_1,y_2)\in\R^3$. Analogous to the planar case, the critical manifold
\begin{equation}
	\cC_0=\left\{ (x,y_1,y_2)\in\R^3\,|\, f(x,y_1,y_2,0)=0\right\}
\end{equation}
is a non-degenerate folded surface if
\begin{equation}\label{eq:def-fold}
	f(0,0,0,0)=0,\; \frac{\partial f}{\partial x}(0,0,0,0)=0,\;\frac{\partial^2 f}{\partial x^2}(0,0,0,0)\neq0,\;\textnormal{D}_yf(0,0,0,0)\neq0\in\R^{2}.
\end{equation}
Without loss of generality we may assume that $\frac{\partial f}{\partial y_1}(0,0,0,0)\neq0$.  Note that there is a fold line along $\ell_f=\left\{(x,y)\in\R^3\,|\, x=0, \, y_1=0\right\}$ and that the phase-space of the slow subsystem is $2$-dimensional. Similar to the planar case, the corresponding critical manifold is a folded surface as shown in Figure \ref{fig:desing-example}.

Analogous to the planar case, we say that the origin is a generic fold point if $g_1(0,0,0,0)\neq0$. In this case \eqref{eq:sfs3d} can be rewritten as
\begin{equation}\label{eq:sfs3dfold0}
	\begin{split}
		x'  &= y_1+x^2 + \cO(\ve x,\ve y_1,\ve y_2,\ve^2,x^3,xy_1y_2,x^3y_1),\\
		y_1' &= \ve( 1 + \cO(x,y_1,y_2,\ve) ),\\
		y_2' &= \ve g_2(x,y_1,y_2,\ve).
	\end{split}
\end{equation}

The analysis of \eqref{eq:sfs3dfold0} as its flow passes near $\ell_f$ is similar to that of the generic fold in Section \ref{sec:genericfold}, and the details can be found, for example, in~\cite[Theorem 1]{szmolyan2004relaxation}. On the other hand, we say that the origin of \eqref{eq:sfs3d} is a \emph{folded singularity} if besides \eqref{eq:def-fold} we have
\begin{equation}
	\left(\frac{\partial f}{\partial y_1}g_1 + \frac{\partial f}{\partial y_2}g_2\right)|_{(0,0,0,0)} = 0.
\end{equation}

If the origin is a folded singularity, then \eqref{eq:sfs3d} can be rewritten as~\cite[Proposition 2.1]{szmolyan2001canards}
\begin{equation}\label{eq:sfs3dfold}
	\begin{split}
		x'  &= y_1+x^2 + \cO(\ve x,\ve y_1,\ve y_2,\ve^2,x^3,xy_1y_2,x^3y_1),\\
		y_1' &= \ve(ax+by_2 + \cO(y_1,\ve,y_2^2,xy_2,x^2)),\\
		y_2' &= \ve ( c+\cO(x,y_1,y_2,\ve) ).
	\end{split}
\end{equation}

Since a folded singularity is necessary for the existence of canards, we now focus on \eqref{eq:sfs3dfold}. Note that due to the leading order terms of \eqref{eq:sfs3dfold} we have that $\cC_0$ is attracting for $x<0$ and repelling for $x>0$. The relationship between the flow on the critical manifold and the fold line is of vital importance. A classification of the possible phase-portraits on $\cC_0$ near $\ell_f$ can be found in~\cite[Theorem 5.1]{takens1976constrained} and in~\cite[Lemma 2.1]{szmolyan2001canards}. Such a classification makes use of the so-called \emph{ desingularized system}~\cite[Section 4.10]{takens1976constrained} which, up to leading order terms, is given by
\begin{equation}\label{eq:desing}
	\begin{split}
		\dot x &= ax+by_2,\\
		\dot y_2 &= -2cx.
	\end{split}
\end{equation}

Note that the origin is an equilibrium point of \eqref{eq:desing}. Even though the origin is not an equilibrium point of \eqref{eq:sfs3dfold}, it is possible to relate the flow of \eqref{eq:desing} to the flow of \eqref{eq:sfs3dfold} on $\cC_0$.  Thus, a classification of the flow of \eqref{eq:desing} is useful and is summarized in Table~\ref{table:class}.

\begin{table}[htbp]\centering
	\begin{tabular}{l l l}
		Condition & Eigenvalues & Name\\
			\hline
		$bc<0$ & $\lambda_{1,2}\in\R$, $\lambda_1<0<\lambda_2$ & saddle\\
		$0<8bc< a^2$ & $\lambda_{1,2}\in\R$, $\lambda_1\lambda_2>0$  & node\\
		$a^2<8bc$& $\lambda_{1,2}\in\mathbb C$, $\Re(\lambda_{1,2})\neq0$ & focus\\
		$bc = 0$, $c\neq0$ & $\lambda_1=0$, $\lambda_2=a$ & saddle-node type I\\
		$bc = 0$, $b\neq0$ & $\lambda_1=0$, $\lambda_2=a$ & saddle-node type II\\
		$8bc = a^2$ & $\lambda_{1,2}=\frac{a}{2}$ & degenerate node\\
	\end{tabular}
	\caption{Topological Classification of \eqref{eq:desing}.}
	\label{table:class}
\end{table}

We can see from the above classification that the saddle, node, and focus cases are generic, while the saddle-node and the degenerate node are both degenerate and of codimension $1$. From \eqref{eq:desing} one can obtain the flow on $\cC_0$ as follows: away from $\ell_f$, the trajectories on $\cC_0|_{x<0}$ are given by trajectories of \eqref{eq:desing} restricted to $x<0$, while  trajectories on $\cC_0|_{x>0}$ are given by reversing the trajectories of \eqref{eq:desing} restricted to $x>0$, see the details in e.g. \cite[Section 4.10]{takens1976constrained}, and Figure \ref{fig:desing-example} for a couple of examples.

Due to the above classification and relationship between the slow subsystem and the desingularized system, we call the origin of \eqref{eq:sfs3dfold} a \emph{folded saddle}, \emph{folded node}, \emph{folded focus}, \emph{folded saddle-node}, or \emph{folded degenerate node} according to the coefficients of Table~\ref{table:class}. We now have the following definition of canards coming from the flow on $\cC_0$.

\begin{definition}
Solutions of the slow subsystem passing through a canard point from an attracting critical manifold to a repelling critical manifold are called \emph{singular canards}. Solutions of the slow subsystem passing through a canard point from a repelling critical manifold to an attracting critical manifold are called \emph{faux canards}.
\end{definition}

See, for example Figure \ref{fig:desing-example} (a) where $\Gamma_1$ is a singular canard and $\Gamma_2$ a faux canard.

It follows from the classification of Table \ref{table:class} and the resulting reduced flow that there are singular canards for all cases, except for the folded focus case. We can now recall the main result of \cite{szmolyan2001canards}, stating under which conditions singular canards persist for $\ve>0$ sufficiently small.

\begin{theorem}[{\cite[Theorem 4.1]{szmolyan2001canards}}]\label{thm:3dcanard1} Assume \eqref{eq:sfs3dfold}. In the folded saddle and in the folded node case singular canards $\Gamma_1$ perturb to maximal canard solutions for sufficiently small $\ve$.

For a folded node with $\lambda_1<\lambda_2<0$ a maximal canard solution corresponding to the weak eigendirection ($\Gamma_2$) exists for sufficiently small $\ve$ provided that the ratio $\mu_2=\lambda_1/\lambda_2$ is not a natural number.

\end{theorem}

See Figure \ref{fig:desing-example} for a schematic representation of the trajectories $\Gamma_1$ and $\Gamma_2$ in each case.  We shall denote the perturbation of $\Gamma_1$ and of $\Gamma_2$ by $\Gamma_{1,\ve}$ and $\Gamma_{2,\ve}$,  respectively. The idea of the proof is similar to that in Section \ref{sec:planarcanards}. The relevant charts are $K_1=\left\{ \bar y_1= 1 \right\}$ and $K_2=\left\{ \bar \ve=1 \right\}$. It is precisely in chart $K_2$ where transversality of the intersection of the slow manifolds can be shown. We now sketch the main procedure: first of all, the blown up vector field in chart $K_2$ has (up to appropriate conditions met by the hypothesis of Theorem \ref{thm:3dcanard1}) explicit algebraic solutions $\gamma_1$, $\gamma_2$ \cite[Lemma 4.2]{szmolyan2001canards}. These are crucial as they connect the repelling and attracting parts of the slow manifolds (equivalently, and analogously to Section \ref{sec:planarcanards}, of centre manifolds in chart $K_1$). Moreover, using variational arguments, it is further shown by studying a Weber equation (see ~\cite[Chapter 19]{abramowitz1964handbook} and \cite[Chapter 12]{olver2010nist}) that the invariant manifolds in chart $K_2$ corresponding to the slow manifolds $\cS_\ve^\txta$ and $\cS_\ve^\txtr$ intersect transversaly along the $\gamma_i$ solutions. Furthermore. the meaning of the ratio $\mu_2$ is quite interesting.

\begin{lemma}[{\cite[Lemma 4.4]{szmolyan2001canards}}] Suppose we are in the node scenario and that $n-1<\mu_2<n$. Then the slow manifolds $\cS_\ve^\txta$ and $\cS_\ve^\txtr$ twist $n-1$ times around the corresponding maximal canard in the neighborhood of the fold curve.
\end{lemma}

\ch{ 
Regarding the folded node case, it is further shown in~\cite{wechselberger2005existence} that one may find more canards bifurcating, upon variation of $\mu_2$, from the weak maximal canard. In terms of the slow manifolds, this means that $\cS_\ve^\txta$ and $\cS_\ve^\txtr$ have further transversal intersections, called \emph{secondary maximal canards}. To better grasp what happens in this situation, let us first note that in the folded node case we can differentiate two types of singular canards: a) singular canards $\Gamma_1$ and $\Gamma_2$ corresponding to the (strong and weak) eigendirections of \eqref{eq:desing} and b) all other singular canards within the shaded region of Figure \ref{fig:desing-example}. The latter type of singular canards may also correspond to maximal canards for $\ve>0$ sufficiently small, this motivates the next definition.

\begin{definition} All maximal canards not obtained as perturbations of the singular canards $\Gamma_1$, $\Gamma_2$ using Theorem \ref{thm:3dcanard1} are called \emph{secondary maximal canards}.
\end{definition}

Skipping the technicalities, the main result of \cite{wechselberger2005existence} is as follows: ``If $\mu_2\in\mathbb N$, then for $\ve>0$ sufficiently small there are additional secondary maximal canards bifurcating from the weak maximal canard upon variation of the ratio $\mu_2$". The proof of the aforementioned bifurcation result requires a highly non-trivial analysis partially contained in \cite{wechselberger2005existence} and further completed in \cite{mitry2017folded}.

More recently, the folded saddle case has been studied in more detail in \cite{mitry2017folded}. Interestingly, it is proven that for certain values of the ratio $\mu_2$ there is rotational behaviour of trajectories around the faux canard $\Gamma_2$ \cite[Lemma 1]{mitry2017folded}. This rotational behaviour has important consequences, in particular, it is shown that trajectories with distinct rotational behaviour can be identified. The boundaries of sectors of trajectories with the same rotational behaviour are given by the so-called \emph{fast manifolds} (the fast manifolds are defined as the nonlinear stable and unstable fibre bundles of the primary faux canard $\Gamma_2$). Furthermore, the transversal intersection of these fast manifolds is responsible for the creation of secondary faux canards. This contrasts with the folded node case where, as described above, the intersection of the slow manifolds of the weak maximal canard is the responsible for the creation of secondary canards. The main analysis of the secondary canards and their bifurcations relies on a variational method and an extension of Melknikov's method.

}

Visualizing the slow manifolds, especially as they twist in the folded node case, is difficult and the reader is referred to e.g.\cite{desroches2008geometry,guckenheimer2005canards} for further details. We finalize this section by briefly sketching some generalizations and referring to works regarding the degenerate saddle-node cases.
\begin{itemize}[leftmargin=*]
	\item The canard theory sketched above can be extended to fast-slow systems with $m\geq2$ slow and $n\geq1$ fast variables \cite{wechselberger2012propos}. The arguments use the assumption of $n-1$ stable fast directions and centre manifold reduction.
	\item The codimension $1$ folded saddle node cases are studied with the help of the blow-up method in \cite{krupa2010local} (FSN-II) and in \cite{vo2015canards} (FSN-I). One primary interest in studying such folded singularity is their relation to complex oscillatory motion such as Mixed Mode Oscillations (MMOs)~\cite{desroches2012mixed} and chaos. Furthermore, folded saddle-nodes appear in many models of applied interest such as the forced van der Pol oscillator~\cite{bold2003forced,guckenheimer2006chaotic,burke2016canards}, biochemical reactions~\cite{milik1998geometry,milik2001multiple,moehlis2002canards}, and neuron dynamics~\cite{ermentrout2009canards,krupa2008mixed,de2015neural}, to mention a few. Unfolding such degenerate cases, that is considering one-parameter families of fast-slow systems with folded singularities, provides a much richer geometry. One difference to keep in mind between the two cases is that the FSN-II has an associated singular Hopf bifurcation with two slow variables~\cite{guckenheimer2008singular}, while the FSN-I does not present a Hopf bifurcation but rather a true saddle-node. 

\end{itemize}

\subsection{Delayed Hopf bifurcation}\label{sec:delayedHopf}

In dynamic bifurcation problems, that is where the bifurcation parameter is slowly changed, it can be observed that the transition to instability occurs for parameter values way beyond the prediction of the static bifurcation diagram \cite{Neishtadt2009}. This effect is known as \emph{delayed loss of stability}. One landmark example of such a phenomenon occurs in the so-called \emph{slow passage through a Hopf bifurcation}~\cite{Baer1989slow}. The aforementioned dynamical behavior is also highly relevant in many applications, e.g.~\cite{strizhak1996slow,barreto2011ion}; see also \cite{kaper2018delayed} for a recent example of slow passage through a Hopf bifurcation in reaction-diffusion equations.

Although the problem of slow passage through a Hopf bifurcation has been long studied, \cite{HAYES20161184} presents a new geometric proof based on the blow-up method, which we now sketch. Consider a fast-slow system \eqref{eq:sfs1}-\eqref{eq:sfs2} with $(m,n)=(2,1)$ given by
	\begin{equation}\label{eq:delayedhopf}
		\begin{split}
			x_1' &= yx_1 - x_2 - x_1(x_1^2+x_2^2),\\
			x_2' &= x_1 + yx_2 - x_2(x_1^2+x_2^2),\\
			y' &= \ve.
		\end{split}
	\end{equation}

Note that for $\ve=0$ we have that \eqref{eq:delayedhopf} corresponds to the normal form of a (supercritical) Hopf bifurcation with $y$ the bifurcation parameter and $y=0$ the bifurcation point. This motivates the following definition.

\begin{definition}
A Hopf bifurcation that occurs in the layer equation of a fast-slow system in which the slow variable acts as the bifurcation parameter is called \emph{Delayed Hopf Bifurcation}.
\end{definition}

The idea now is to study the dynamics of  \eqref{eq:delayedhopf} as the ``dynamic parameter'' $y$ slowly passes through the Hopf bifurcation point $y=0$. It is straightforward to see that the critical manifold is
\begin{equation}
	\cC_0 = \left\{ (x_1,x_2,y)\in\R^3\, | \, x_1=0,x_2=0 \right\},
\end{equation}
i.e., it is just the $y$-axis. One can further show, for example by changing to polar coordinates, that $\cC_0|_{y<0}$ is attracting while $\cC_0|_{y>0}$ is repelling, and that the origin is a non-hyperbolic point of the layer equation. The next observation is of crucial importance: the linearization of the layer equation \eqref{eq:delayedhopf} at the non-hyperbolic origin is not nilpotent. In fact it is given by $\dot x =Jx$, where
\begin{equation}
		J=\begin{bmatrix}
			0 & -1\\
			1 & 0
	\end{bmatrix}.
\end{equation}
This is a first obstacle to use the blow-up method. Furthermore, note that such linearization corresponds to a fast rotation of the fast variables. We will see below that the main idea of \cite{HAYES20161184} is to ``eliminate'' such a fast rotation allowing the use of the blow-up method. It is convenient to consider a modified complex version
\begin{equation}\label{eq:delayedhopf2}
	\begin{split}
		z' &= (y+\imath)z - |z|^2z + \ve k_0,\\
		y' &= \ve.
	\end{split}
\end{equation}
System \eqref{eq:delayedhopf2}, with $k_0=0$ is obtained from \eqref{eq:delayedhopf} via the change of coordinates $z=x_1+\imath x_2$ with $\imath=\sqrt{-1}$. The term $\ve k_0$, $k_0\neq0$, is added to break the invariance\footnote{In principle a generic term $\cO(\ve)$ can be added instead of $\ve k_0$, but the results are similar in such a case.} of $\cC_0=\left\{ z =0 \right\}$ for $\ve>0$. It is also clear from \eqref{eq:delayedhopf2} that $z=0$ is a stable (resp. unstable) focus for $y<0$ (resp. $y>0$) and that at $z=0$ an exchange of stability occurs. Let $\delta>0$ be small and, as usual, define
\begin{equation}
	\begin{split}
		\cS_0^\txta &= \left\{ (z,y)\in\cC_0\,| \, y<-\delta \right\},\\
		\cS_0^\txtr &= \left\{ (z,y)\in\cC_0\,| \, y>\delta \right\},
	\end{split}
\end{equation}
and we denote by $\cS_\ve^\txta$ and $\cS_\ve^\txtr$ the corresponding slow manifolds. The manifolds $\cS_0^\txta$ and $\cS_0^\txtr$ have analytic extensions to $y=0$, and the goal is to compute the distance between $\cS_\ve^\txta$ and $\cS_\ve^\txtr$ at $y=0$ for $\ve>0$ sufficiently small. Let such a distance be denoted by $d(\ve)$. The main result of \cite{HAYES20161184} is a new geometric proof of the following well-known theorem.

\begin{theorem}[{\cite[Theorem 3.1]{HAYES20161184} see also \cite{neishtadt1987persistence,shishkova1973analysis}}]\label{thm:hayes} Let $k_0\neq 0$. For sufficiently small values of $\ve>0$,
	\begin{equation}
		d(\ve)\leq \left( \sqrt\ve\sqrt{2\pi}\kappa+\cO(\ve) \right)\exp\left( -\frac{1}{2\ve}\right),
	\end{equation}
where $\kappa>0$ is bounded away from $0$, and is determined in part by $|k_0|$.
\end{theorem}

\begin{figure}[htbp]\centering
	\begin{tikzpicture}
		\node at (0,0){
		\includegraphics[scale=1.25]{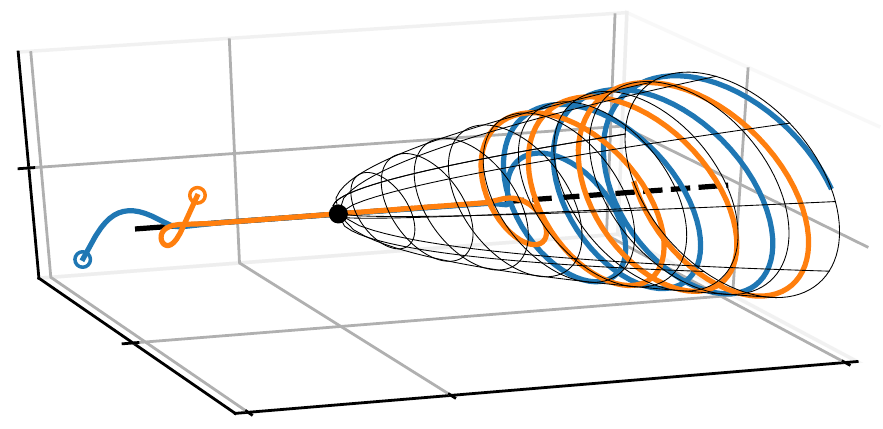}
		};
		\node at (-1.6,.25) {$H$};
		\node at (0.2,-2.6) {$0$};
		\node at (-2.2,-2.9) {$-1$};
		\node at ( 5.3,-2.2) {$2$};
		\node at (-4.5,-1.75) {$0$};
		\node at (-5.75,.55) {$0$};
		\node at (-6.25,.55) {$x_2$};
		\node at (-5,-1.75) {$x_1$};
		\node at (1,-3) {$y$};
	\end{tikzpicture}
	\caption{Simulation of a slow passage through a Hopf bifurcation of \eqref{eq:delayedhopf} with $\ve=0.05$. The critical manifold is given by $\cC_0=\left\{ (x_1,x_2)=(0,0) \right\}$, it is attracting (solid black line) for $y<0$ and repelling (dashed black line) for $y>0$. The point $H=(0,0,0)$ denotes the bifurcation point and the wire-frame paraboloid represents the amplitude of oscillations according to the $y$-parameter value. We show two trajectories and observe that the onset of oscillations occurs way beyond the Hopf point.}
	\label{fig:hopf}
\end{figure}
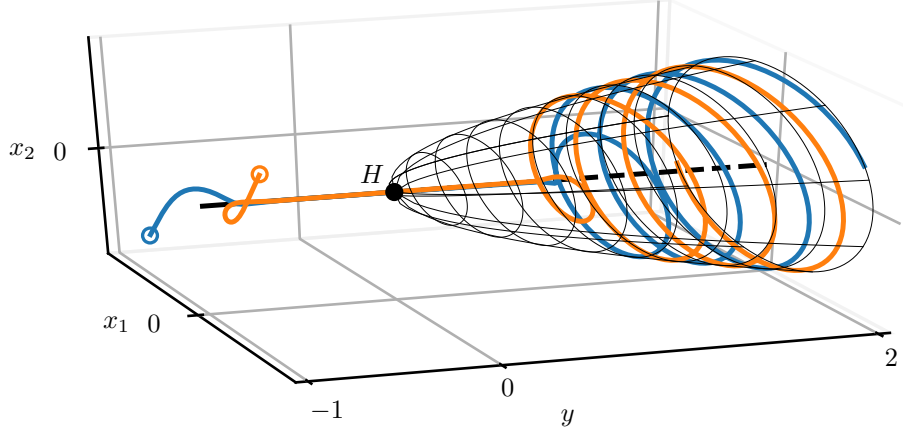

The implication of Theorem \ref{thm:hayes} is that, due to the fact that the separation between the slow manifolds $\cS_\ve^\txta$ and $\cS_\ve^\txtr$ is exponentially small, solutions starting near $\cS_\ve^\txta$ stay close to $\cS_\ve^\txtr$ for time $\cO(1/\ve)$ after crossing the bifurcation point. This is observed as a delay in the onset of (large amplitude) oscillations. To prove Theorem \ref{thm:hayes} one first considers $y$ in the complex plane, that is $y=y_R+\imath y_I\in\mathbb C$. Next, one chooses a particular path, namely
\begin{equation}
	\Gamma_{-1}:y = s-\imath, \qquad s\in(-\infty,\infty).
\end{equation}
The main observation is that, since \eqref{eq:delayedhopf2} is analytic, by Cauchy's integral theorem it is the equivalent to integrate along $\Gamma_{-1}$ than to integrate along the real $y$-axis. The advantage, however, is that along $\Gamma_{-1}$, the fast rotation is eliminated, and the resulting non-hyperbolic equilibrium point is nilpotent. So, the fast-slow system \eqref{eq:delayedhopf2} along $\Gamma_{-1}$ reads as
\begin{equation}\label{eq:delayedhopf3}
	\begin{split}
		z' &= sz-|z|^2z+\ve k_0,\\
		s' &= \ve,
	\end{split}
\end{equation}
where we can see that $(z,s,\ve)=(0,0,0)$ is indeed a nilpotent point. We now proceed to sketch the blow-up analysis. The blow-up map is given by
\begin{equation}
	z = r\bar z, \quad s = r\bar s, \quad \ve = r^2\bar\ve.
\end{equation}
The relevant charts are $K_1 = \left\{ \bar s=-1\right\}$, $K_2 = \left\{ \bar\ve=1 \right\}$, and $K_3 = \left\{ \bar s=1\right\}$. Next, the strategy is similar to the planar fold case: in chart $K_1$ (resp. $K_3$) one finds a $2$-dimensional centre-stable (resp.~centre-unstable) manifold corresponding to $\cS_\ve^\txta$ (resp. $\cS_\ve^\txtr$). Then, in chart $K_2$ one tracks such invariant manifolds. To compute the distance between the manifolds we take advantage of the fact that the blown-up vector field in $K_2$ can be integrated along appropriate paths. We remark that due to the cubic terms in \eqref{eq:delayedhopf3} a secondary blow-up is performed. The delayed onset of oscillations due to a slow passage through a Hopf bifurcation can be observed in Figure \ref{fig:hopf}.

\subsection{Beyond fold and Hopf singularities}\label{sec:beyondfold}

As we have reviewed so far, the understanding of codimension-one fold and Hopf points is of extreme importance due to their genericity in nonlinear fast-slow systems. However, there are many other singularities that have also been studied via the blow-up method, which we summarize here.

\begin{description}[leftmargin=0cm]
	\item[Transcritical and pitchfork singularities] Let us consider \eqref{eq:sfs2} with $(x,y)\in\R^2$ and assume that the layer equation undergoes either a transcritical or a pitchfork bifurcation at the origin where $y\in\R$ is the bifurcation parameter~\cite{krupa2001extendingtrans}. The goal is to understand the dynamics for $\ve>0$ small in the neighbourhood of the aforementioned singularities.

\emph{In the case of the transcritical singularity}, it can be shown (Lemma 2.1 in~\cite{krupa2001extendingtrans}) that a canonical form is given as
\begin{equation}\label{eq:trans}
		\begin{split}
			x' &= x^2 - y^2 + \omega\ve + \cO(x^3,x^2y,xy^2,y^3,\ve x, \ve y, \ve^2),\\
			y' &= \ve(1+\cO(x,y,\ve)),
		\end{split}
\end{equation}
where $\omega$ is considered as a continuous parameter. The role of the parameter $\omega$ is to distinguish two types of dynamics, as is stated below in Theorem \ref{thm:trans}. The critical manifold $\cC_0$ is given by $\cC_0=\left\{ (x,y)\in\R^2 \, | \, x=\pm y + \cO(y^2) \right\}$, which clearly has a self-intersection at the origin. From another perspective, the critical manifold can be regarded as the union of one dimensional manifolds $\left\{ x=y+\cO(y^2)\right\}$  and $\left\{ x=-y+\cO(y^2)\right\}$, which intersect transversally at the origin. Such intersection implies that the origin is a non-hyperbolic equilibrium point of the layer equation. In either case, it is convenient to assume that a compact subset of $\cC_0$ around the origin is given by the union of four branches, namely $\cS_{0,+}^\txta,\,\cS_{0,-}^\txta,\,\cS_{0,+}^\txtr,\,\cS_{0,-}^\txtr$, as shown in Figure \ref{fig:trans}. Thus, away from the origin, Fenichel's theorem guarantees that the previous critical manifolds persist as slow manifolds $\cS_{\ve,+}^\txta,\,\cS_{\ve,-}^\txta,\,\cS_{\ve,+}^\txtr,\,\cS_{\ve,-}^\txtr$ for $\ve>0$ sufficiently small. The idea of~\cite{krupa2001extendingtrans} is to describe the fate of $\cS_{\ve,-}^\txta$ as it crosses a small neighborhood of the origin. In order to state the main result regarding the transcritical singularity, let $\rho>0$ be small and $I\in\R$ a small interval containing $0\in\R$. Define the sections
	\begin{equation}
		\begin{split}
			\Sigma^{\EN} &= \left\{ (x,y)\in\R^2\, | \, x=-\rho, \, y+\rho\in I \right\},\\
			\Sigma_e^{\EX} &= \left\{ (x,y)\in\R^2\, | \, x=\rho, \, y\in I \right\},\\
			\Sigma_a^{\EX} &= \left\{ (x,y)\in\R^2\, | \, x=-\rho, \, y-\rho\in I \right\}.
		\end{split}
	\end{equation}
Let $\Pi_e:\Sigma^{\EX}\to\Sigma_e^{\EX} $ and $\Pi_a:\Sigma^{\EN}\to\Sigma_a^{\EX} $ denote the corresponding transitions.

\begin{theorem}[\cite{krupa2001extendingtrans}]\label{thm:trans} Consider \eqref{eq:trans} and fix $\omega\neq1$. There exists $\ve_0>0$ such that the following assertions hold for $\ve\in(0,\ve_0]$.
	\begin{itemize}
		\item[(T1)] If $\omega>1$ then the manifold $\cS_{\ve,-}^\txta$ passes through $\Sigma_e^{\EX}$ at a point $(\rho,h(\ve))$ where $h(\ve)\in \cO(\sqrt\ve)$. The section $\Sigma^{\EN}$ is mapped by $\Pi_a$ to an interval containing $\cS_{\ve,-}^\txta\cap\Sigma_a^{\EX}$ of size $\cO(\exp(-C/\ve))$, where $C$ is a positive constant.
		\item[(T2)] If $\omega<1$ then $\Sigma^{\EN}$ (including the point $\cS_{\ve,-}^\txta\cap\Sigma^{\EN}$) is mapped by $\Pi_e$ to an interval about $\cS_{\ve,+}^\txta$ of size $\cO(\exp(-C/\ve))$, where $C$ is a positive constant.
	\end{itemize}
	\end{theorem}

See Figure \ref{fig:trans} for a schematic of the results of Theorem \ref{thm:trans}. It is important to note that the constant $\ve_0$ tends to $0$ as $\omega\to 1$. In fact, the case $\omega=1$ is degenerate and is not detailed in \cite{krupa2001extendingtrans}. However, it is argued in Remark 2.2 of \cite{krupa2001extendingtrans} that canard solutions occur in this case, the treatment being similar to that of \cite{krupa2001extending}, see also the Example in Section 3 of \cite{krupa2001extendingtrans}. From another perspective, i.e. not using blow-up, a particular case of this degenerate scenario is investigated in \cite{schecter1985persistent}.

For the proof of Theorem \ref{thm:trans}, one employs a blow-up map defined by $\phi(r,\bar x, \bar y, \bar\ve)=(r\bar x, r\bar y, r^2\bar\ve)$. The blow-up analysis is carried out in three charts, namely $K_1=\left\{ \bar x = -1\right\}$, $K_2=\left\{ \bar \ve = 1\right\}$ and $K_3=\left\{ \bar x = 1\right\}$. The analysis in charts $K_1$ and $K_3$ is similar and amounts to describing the organization of the dynamics according to two $2$-dimensional centre manifolds, and which are tracked along the blow-up space. Another important remark is that in the charts $K_1$ and $K_3$ the parameter $\omega$ does not play any essential role. In chart $K_2$, the main object of study is the ODE
\begin{equation}
	\begin{split}
		x_2' &= x^2 - y^2 + \omega,\\
		y_2' &= 1,
	\end{split}
\end{equation}
where as usual we are using $(x_2,y_2)$ to denote local coordinates. The main goal in this chart is to connect the appropriate centre manifolds of charts $K_1$ and $K_3$.

A more general scenario compared to the one presented above is treated in \cite{de2015planar}, where in brief terms, unfoldings of \eqref{eq:trans} are considered. In fact, since \cite{de2015planar} considers fast-slow systems with a transcritical singularity in its most generic context, more passages as compared to those shown in Figure \ref{fig:trans} can occur,  \cite[Theorems 5 and 6]{de2015planar}. The most important of such passages is, of course, the canard. A consequence of this analysis is the extension of known results on the stability of canard cycles~\cite{de2011cyclicity} to the case including finite passages through a transcritical singularity~\cite[Theorem 4]{de2015planar}.\medskip

Regarding the \emph{pitchfork singularity}, it is shown (Lemma 4.1 of \cite{krupa2001extendingtrans}) that a canonical form is given by
\begin{equation}\label{eq:pitch}
		\begin{split}
			x' &= x(y-x^2) + \omega\ve + ax^2y + \cO(xy^2,y^3,\ve x, \ve y, \ve^2),\\
			y' &= \ve(\pm 1 + bx + \cO(y,\ve)),
		\end{split}
\end{equation}
where $\omega$ is assumed non-zero and $(a,b)\in\R^2$. In this case, the analysis and results depend on the sign in \eqref{eq:pitch}. Similar to the transcritical case, the critical manifold can be given as the union of four branches $\cS_{0,\textnormal{t}}^\txta$, $\cS_{0,\textnormal{t}}^\txtr$, $\cS_{0,+}^\txta$, $\cS_{0,-}^\txta$. Away from a small neighbourhood of the origin and for $\ve>0$ sufficiently small, the previous four branches persist as the slow manifolds $\cS_{\ve,\textnormal{t}}^\txta$, $\cS_{\ve,\textnormal{t}}^\txtr$, $\cS_{\ve,+}^\txta$, $\cS_{\ve,-}^\txta$ respectively, see Figure \ref{fig:trans}. Let $\rho>0$ be fixed and $I\subset\R$ be a small interval around $0\in\R$. Define the sections
\begin{equation}
		\begin{split}
			\Sigma^{\textnormal{t}} &=  \left\{ (x,y)\in\R^2\, | \, x\in I, \, y=-\rho \right\},\\
			\Sigma^+ &= \left\{ (x,y)\in\R^2\, | \, x=\rho, \, y-\rho^2\in I \right\},\\
			\Sigma^- &= \left\{ (x,y)\in\R^2\, | \, x=-\rho, \, y-\rho^2\in I \right\},
		\end{split}
\end{equation}
and define transitions $\Pi^t:\Sigma^{\textnormal{t}}\to\Sigma^+\cup\Sigma^-$ and $\Pi^\pm:\Sigma^\pm\to\Sigma^{\textnormal{t}}$. The former transition is well defined when the sign in \eqref{eq:pitch} is positive while the latter is well defined when the sign in \eqref{eq:pitch} is negative. Regarding the pitchfork singularity, the main result is:

\begin{theorem}[\cite{krupa2001extendingtrans}]\label{thm:pitch} Consider \eqref{eq:pitch} and fix $\omega\neq0$. There exists $\ve_0>0$ such that the following assertions hold for $\ve\in(0,\ve_0]$.
		\begin{itemize}
			\item[(P1)] If the sign in \eqref{eq:pitch} is positive and $\omega>0$ then $\Sigma^{\textnormal{t}}$ (including the point $\Sigma^{\textnormal{t}}\cap\cS_{\ve,\textnormal{t}}^\txta$) is mapped by $\Pi^t$ to an interval about $\Sigma^+\cap\cS_{\ve.+}^\txta$ of size $\cO(\exp(-C/\ve))$, where $C$ is a positive constant.
			\item[(P2)] If the sign in \eqref{eq:pitch} is positive and $\omega<0$ then $\Sigma^{\textnormal{t}}$ (including the point $\Sigma^{\textnormal{t}}\cap\cS_{\ve,\textnormal{t}}^\txta$) is mapped by $\Pi^t$ to an interval about $\Sigma^-\cap\cS_{\ve.-}^\txta$ of size $\cO(\exp(-C/\ve))$, where $C$ is a positive constant.
			\item[(P3)] If the sign in \eqref{eq:pitch} is negative then $\Sigma^+$ and $\Sigma^-$ are mapped by $\Pi^+$ and $\Pi^-$, respectively, to intervals about $\cS_{\ve,\textnormal{t}}^\txta\cap\Sigma^{\textnormal{t}}$ of size $\cO(\exp(-C/\ve))$, where $C$ is a positive constant.
		\end{itemize}
\end{theorem}

See Figure \ref{fig:pitch} for a schematic of the \ch{statements} of Theorem \ref{thm:pitch}. In this case, the blow-up map is defined via $\phi(r,\bar x, \bar y, \bar\ve)=(r\bar x, r^2\bar y, r^4\bar\ve)$. The blow-up analysis is carried out in five charts, namely $K_1=\left\{ \bar y = -1\right\}$, $K_2=\left\{ \bar \ve = 1\right\}$, $K_3=\left\{ \bar y = 1\right\}$,  $K_4=\left\{ \bar x = -1\right\}$ and $K_5=\left\{ \bar x = 1\right\}$. As in the transcritical case, the idea is to track centre manifolds along the blow-up space. The relevant centre manifolds are first found in charts $K_1$ and $K_3$ and are related to normally hyperbolic parts of the critical manifold. Accordingly, $\cM_{1,\textnormal{t}}^\txta$, $\cM_{3,+}^\txta$, $\cM_{3,-}^\txta$ and $\cM_{3,\textnormal{t}}^\txtr$ denote centre manifolds (in the blow-up space) associated to the branches $\cS_{0,\textnormal{t}}^\txta$, $\cS_{0.+}^\txta$, $\cS_{0,-}^\txta$ and $\cS_{0,\textnormal{t}}^\txtr$ respectively. On the other hand, in charts $K_4$ and $K_5$ one finds hyperbolic equilibrium points.

The most important analysis occurs in chart $K_2$. In chart $K_2$ one must connect the aforementioned centre manifolds. We sketch the key argument: let $\cM_{2,\textnormal{t}}^\txta$ and $\cM_{2,\textnormal{t}}^\txtr$ denote the centre manifolds $\cM_{1,\textnormal{t}}^\txta$ and $\cM_{3,\textnormal{t}}^\txtr$, respectively, written in the coordinates of chart $K_2$. The main result of the analysis in chart $K_2$ is that the invariant manifolds $\cM_{2,\textnormal{t}}^\txta$ and $\cM_{2,\textnormal{t}}^\txtr$ intersect transversally along $(x_2,y_2,r_2)=(0,0,t_2)$, where $t_2\in(-\infty,\infty)$. The proof of this fact is done via the analysis of a Melknikov integral, see Section 3 in\cite{krupa2001extendingtrans}. In terms of the slow manifolds, the previous result implies that for $\ve>0$ sufficiently small and $\omega\neq 0$, the distance between the slow manifolds $\cS_{\ve,\textnormal{t}}^\txta$ and $\cS_{\ve,\textnormal{t}}^\txtr$, around the origin, is non-zero. On the other hand, when $\omega=0$ the slow manifolds  $\cS_{\ve,\textnormal{t}}^\txta$ and $\cS_{\ve,\textnormal{t}}^\txtr$ are connected.

	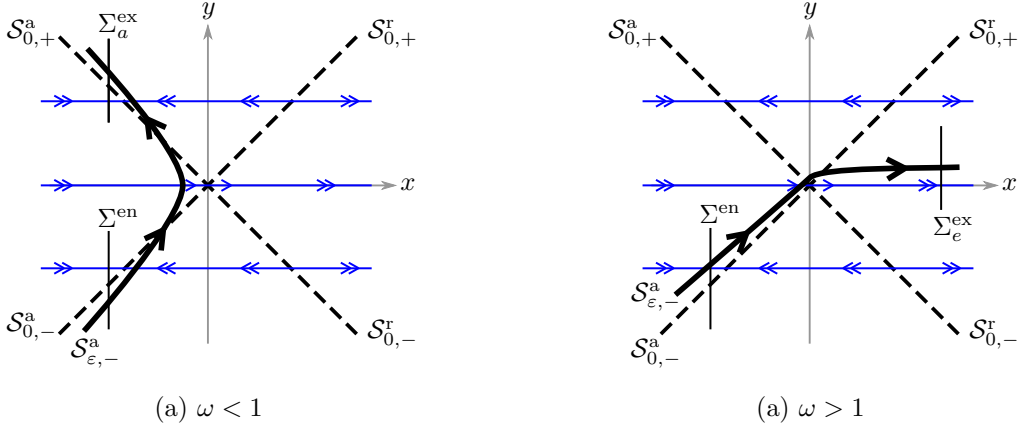
\begin{figure}[htbp]\centering
	\begin{tikzpicture}
		\node at (0,0){
		\includegraphics{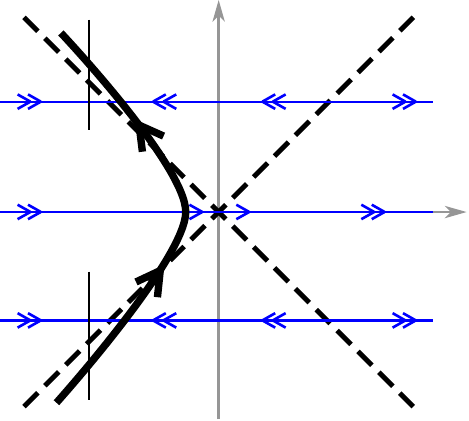}
		};
		\node at (2.5,0) {$x$};
		\node at (-.15,2.3) {$y$};
		\node at (-2.5,-1.9){$\cS_{0,-}^\txta$};
		\node at (-2.5, 2.){$\cS_{0,+}^\txta$};
		\node at ( 2.3,-2.){$\cS_{0,-}^\txtr$};
		\node at ( 2.3, 2.){$\cS_{0,+}^\txtr$};
		\node at ( -1.35, -.45){$\Sigma^\EN$};
		\node at ( -1.35, 2.1){$\Sigma^\EX_a$};
		\node at (-1.65,-2.25){$\cS_{\ve,-}^\txta$};
		\node at (-.15,-3) {(a) $\omega<1$};
		\node[xshift=8cm] at (0,0){
		\includegraphics{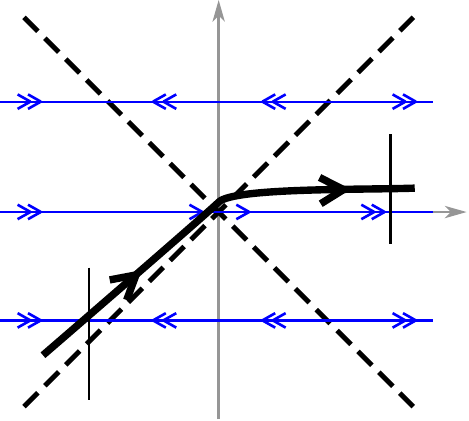}
		};
		\node[xshift=8cm] at (2.5,0) {$x$};
		\node[xshift=8cm] at (-.15,2.3) {$y$};
		\node[xshift=8cm] at (-2.2,-2.25){$\cS_{0,-}^\txta$};
		\node[xshift=8cm] at (-2.5, 2.){$\cS_{0,+}^\txta$};
		\node[xshift=8cm] at ( 2.3,-2.){$\cS_{0,-}^\txtr$};
		\node[xshift=8cm] at ( 2.3, 2.){$\cS_{0,+}^\txtr$};
		\node[xshift=8cm] at ( -1.35, -.4){$\Sigma^\EN$};
		\node[xshift=8cm] at ( 1.75, -.6){$\Sigma^\EX_e$};
		\node[xshift=8cm] at (-2.2,-1.5){$\cS_{\ve,-}^\txta$};
		\node[xshift=8cm] at (-.15,-3) {(a) $\omega>1$};
	\end{tikzpicture}
		\caption{Flow of a fast-slow system near a transcritical singularity for $\omega<1$ on the left and $\omega>1$ on the right. The case $\omega=1$ (not shown) corresponds to a maximal canard in which $\cS_{\ve,-}^\txta$ and $\cS_{\ve,+}^\txtr$ coincide.}
		\label{fig:trans}
	\end{figure}

	\begin{figure}[htbp]\centering
	\begin{tikzpicture}
		\node at (0,0){
		\includegraphics{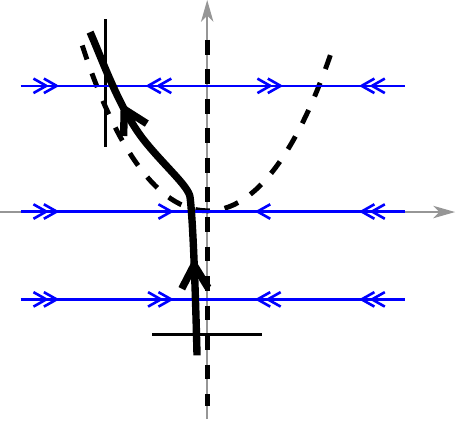}
		};
		\node at (2.5,0) {$x$};
		\node at (-.15,2.3) {$y$};
		\node at ( .55,-1.25) {$\Sigma^{\textnormal{t}}$};
		\node at ( -1.1,2.15) {$\Sigma^-$};
		\node at ( 0.1,-1.95) {$\cS_{0,\textnormal{t}}^\txta$};
		\node at ( 0.1, 1.7) {$\cS_{0,\textnormal{t}}^\txtr$};
		\node at ( 1.4, 1.6) {$\cS_{0,+}^\txta$};
		\node at ( -1.8, 1.6) {$\cS_{0,-}^\txta$};
		\node at ( -.75,-.5) {$\cS_{\ve,\textnormal{t}}^\txta$};
		\node at (-.15,-3) {(a)};
		\node[xshift=8cm] at (0,0){
		\includegraphics{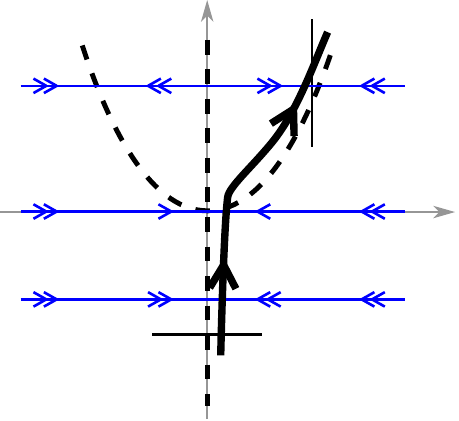}
		};
		\node[xshift=8cm] at (2.5,0) {$x$};
		\node[xshift=8cm] at (-.15,2.3) {$y$};
		\node[xshift=8cm] at ( .55,-1.25) {$\Sigma^{\textnormal{t}}$};
		\node[xshift=8cm] at (  .95,2.15) {$\Sigma^+$};
		\node[xshift=8cm] at ( 0.1,-1.95) {$\cS_{0,\textnormal{t}}^\txta$};
		\node[xshift=8cm] at ( 0.1, 1.7) {$\cS_{0,\textnormal{t}}^\txtr$};
		\node[xshift=8cm] at ( 1.4, 1.6) {$\cS_{0,+}^\txta$};
		\node[xshift=8cm] at ( -1.8, 1.6) {$\cS_{0,-}^\txta$};
		\node[xshift=8cm] at ( .3,-.5) {$\cS_{\ve,\textnormal{t}}^\txta$};
		\node[xshift=8cm] at (-.15,-3) {(b)};
		\node[xshift=4cm,yshift=-6cm] at (0,0){
		\includegraphics{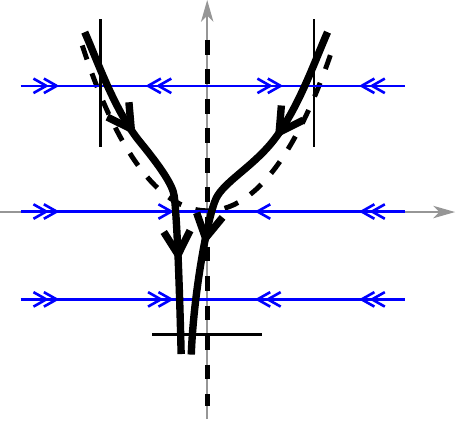}
		};
		\node[xshift=4cm,yshift=-6cm] at (2.5,0) {$x$};
		\node[xshift=4cm,yshift=-6cm] at (-.15,2.3) {$y$};
		\node[xshift=4cm,yshift=-6cm] at ( .55,-1.25) {$\Sigma^{\textnormal{t}}$};
		\node[xshift=4cm,yshift=-6cm] at ( -1.1,2.15) {$\Sigma^-$};
		\node[xshift=4cm,yshift=-6cm] at (  .95,2.15) {$\Sigma^+$};
		\node[xshift=4cm,yshift=-6cm] at ( -.9,-.5) {$\cS_{\ve,-}^\txta$};
		\node[xshift=4cm,yshift=-6cm] at ( .25,-.3) {$\cS_{\ve,+}^\txta$};
		\node[xshift=4cm,yshift=-6cm] at (-.15,-3) {(c)};
	\end{tikzpicture}
		\caption{Phase-portraits of a fast-slow system near a pitchfork singularity for the cases of Theorem \ref{thm:pitch}. Phase-portrait (a) (resp. (b)) corresponds to the case $\omega<0$ (resp.~$\omega>0$) with the sign of \eqref{eq:pitch} positive, while phase-portrait (c) corresponds to the sign of \eqref{eq:pitch} negative.}
		\label{fig:pitch}
	\end{figure}

	$ $\bigskip

\item[Cusp singularity] The cusp singularity is more degenerate than the fold. In fact, at the level of classification of singularities of generic smooth functions, the cusp singularity is the next in the list of singularities (or catastrophes)~\cite{arnold1990theory}. This means that the cusp is of codimension $2$, that is, it appears generically in $2$-parameter families of smooth functions.

In the context of fast-slow systems, the cusp point was first studied in \cite{broer2013geometric}. Given a fast-slow system \eqref{eq:sfs1} with $(m,n)=(1,2)$, the conditions for having a generic cusp point at the origin are
\begin{equation}\label{eq:cusp1}
	\begin{split}
			f(0,0,0,0)=0,\quad \frac{\partial f}{\partial x}(0,0,0,0)=0,\quad
			\frac{\partial^2 f}{\partial x^2}(0,0,0,0)=0
	\end{split}
\end{equation}
and
\begin{equation}\label{eq:cusp2}
\begin{split}
		\frac{\partial^3 f}{\partial x^3}(0,0,0,0)\neq 0 \quad \begin{bmatrix}
			\frac{\partial f}{\partial y_1}(0,0,0,0) & \frac{\partial f}{\partial y_2}(0,0,0,0)
\end{bmatrix}\cdot g(0,0,0,0)\neq0.
\end{split}
\end{equation}
If a fast-slow system satisfies \eqref{eq:cusp1}-\eqref{eq:cusp2}, then it can be written as~\cite[Proposition 2.5]{broer2013geometric}
\begin{equation}\label{eq:cusp3}
	\begin{split}
		x' &= -(x^3 + y_2x+y_1) + \cO(\ve,xy_1,x^3y_2,x^4),\\
		y_1' &= \ve(1+ \cO(x,y_1,y_2,\ve)),\\
		y_2' &= \ve \cO(x,y_1,y_2,\ve).
	\end{split}
\end{equation}
We show in Figure \ref{fig:cusp1} the corresponding critical manifold
\begin{equation}
	\cC_0=\left\{ (x,y_1,y_2)\in\R^3\,|\, -(x^3 + y_2x+y_1) + \cO(\ve,xy_1,x^3y_2,x^4)=0\right\}
\end{equation}
and the reduced flow near the cusp point for a particular choice of signs in \eqref{eq:cusp2}. To state the main theorem of \cite{broer2013geometric} let us define the sections
\begin{equation}
	\begin{split}
		\Sigma^{\EN} &= \left\{ (x,y_1,y_2)\in\R^3\, | \, y_1=-\delta \right\},\\
		\Sigma^{\EX} &= \left\{ (x,y_1,y_2)\in\R^3\, | \, y_1=\delta \right\},
	\end{split}
\end{equation}
where $\delta>0$ is small. Furthermore, let $\cS_0^-$ and $\cS_0^+$ be compact subsets of $\cC_0$ satisfying\footnote{Note that both critical manifolds $\cS_{0,-}$ and $\cS_{0,+}$ are attracting, so we omit the superscript ``$\txta$''.}
\begin{equation}
	\begin{split}
		\cS_{0,-} &= \left\{ (x,y_1,y_2)\in\cC_0\,|\, y_1< -\delta \right\},\\
		\cS_{0,+} &= \left\{ (x,y_1,y_2)\in\cC_0\,|\, y_1>\delta \right\}.
	\end{split}
\end{equation}
Accordingly we denote by $\cS_{\ve,-}$ and $\cS_{\ve,+}$ the respective slow manifolds. The main result describes the passage of trajectories starting near the slow manifold $\cS_{\ve,-}$ through a neighborhood of the cusp point:

\begin{theorem}[{\cite[Theorem 2]{broer2013geometric}}]\label{thm:cusp1} For system \eqref{eq:cusp3} there exists $\ve_0>0$ sufficiently small such that for all $0<\ve\leq\ve_0$ the following statements hold:
\begin{itemize}
		\item[(C1)] The transition map $\Pi:\Sigma^{\EN}\to\Sigma^{\EX}$ induced by the flow of \eqref{eq:cusp3} is a diffeomorphism mapping a rectangular neighborhood of $\cS_{0,-}$ into $\Sigma^{out}$.
		\item[(C2)] The choice of $\cS_{\ve,+}$ can be made in such a way that $\Pi(\Sigma^{\EN}\cap\cS_{\ve,-})\subset \cS_{\ve,+}$.
		\item[(C3)] The map $\Pi$ is exponentially contracting in the $x$ direction and the derivative in the $y_2$ direction is bounded above and below. More precisely
		\begin{equation}
			\left| \frac{\partial \Pi}{\partial x} \right|\in \cO(\exp(-C/\ve))
\end{equation}
for some positive constant $C$, and there exists a constant $M>0$ such that
\begin{equation}
	\frac{1}{M}\leq\left| \frac{\partial \Pi}{\partial y_2}\right|\leq M.
\end{equation}
\end{itemize}
\end{theorem}

The proof of Theorem \ref{thm:cusp1} follows the same ideas as for the generic fold case, so we point out some of the differences: first of all, due to the quasihomogeneity of the fast equation, the blow-up map is given by $(x,y_1,y_2,\ve)=(\br^3\bx,\br^2\by_1, \br\by_2,\br^5\be)$. These weights can be obtained by checking the quasi-homogeneity of the vector field. Next, due to the geometric properties of the critical manifold, it is now necessary to study the blown-up vector field in five charts, namely: $K_{\pm\by_1}=\left\{ \by_1=\pm 1 \right\}$, $K_{\pm\by_2}=\left\{ \by_2=\pm 1 \right\}$ and $K_{\be}=\left\{ \be=1 \right\}$. The analysis starts in $K_{-\bar y_1}$. Then, depending on the initial conditions one transitions to either $K_{+\bar y_2}$, $K_{ \bar \ve}$, or to $K_{-\bar y_2}$. In $K_{+\bar y_2}$ we study trajectories passing sufficiently away from the cusp point and along the regular part of $\cS_0$. In $K_{\bar \ve}$ we study trajectories passing through a small neighborhood of the cusp point. In $K_{-\bar y_2}$ we study trajectories passing sufficiently away from the cusp point and across the folded part of $\cS_0$.  It is important to remark that the flow is regular in the charts $K_{+\bar y_2}$, $K_{\bar \ve}$ and $K_{-\bar y_2}$. Moreover, in $K_{-\bar y_2}$ one may invoke the results of the generic fold in Section \ref{sec:genericfold}. Finally, we further transition from $K_{+\bar y_2}$, $K_{+\bar \ve}$ and $K_{-\bar y_2}$ to $K_{+\bar y_1}$ completing the analysis. A schematic of Theorem \ref{thm:cusp1} is shown in Figure \ref{fig:cusp1}.

Another contribution of \cite{broer2013geometric} is the study of trajectories with initial conditions below the cusp point. Interestingly, it can be shown that, under appropriate assumptions, two trajectories with exponentially close initial conditions cross $\Sigma^{\EX}$ algebraically apart, see \cite[Proposition 3]{broer2013geometric} and Figure \ref{fig:cusp1} for a sketch.

\begin{figure}[htbp]\centering
	\begin{tikzpicture}
		\node at (0,0){
		\includegraphics[scale=1.25]{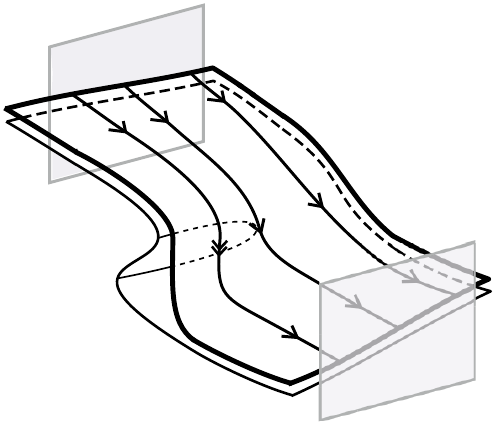}
		};
		\node at (-.5,2.9) {$\Sigma^\EN$};
		\node at (2.75,-2.45) {$\Sigma^\EX$};
		\node at (-3.15,.9) {$\cS_0$};
		\node at (-3.25,1.5) {$\cS_\ve$};
		\fill[black] (0.1,-.2) circle (.075);
		\draw[gray,->] (-2.75,-2.) --++(.5,-.2) node[right,yshift=-.1cm]{\small$y_1$};
		\draw[gray,->] (-2.75,-2.) --++(.5,.2) node[right,yshift=.05cm]{\small$y_2$};
		\draw[gray,->] (-2.75,-2.) --++(0,.5) node[above]{\small$x$};

	\end{tikzpicture}\hfill
	\begin{tikzpicture}
		\node at (0,0){
		\includegraphics[scale=1.25]{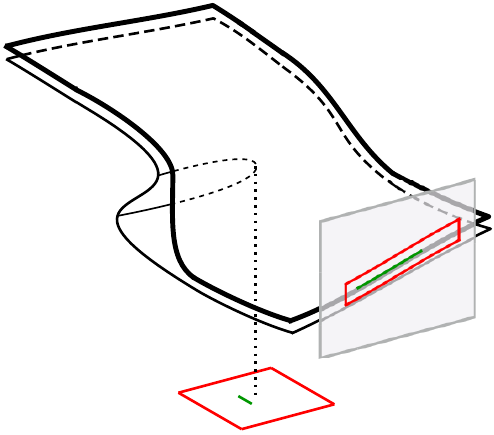}
		};
		\node at (2.75,-2.45 +.85) {$\Sigma^\EX$};
		\node at (-3.15,.9   +.85) {$\cS_0$};
		\node at (-3.25,1.5  +.85) {$\cS_\ve$};
		\node at (1.5,-3.45 +.85) {$\Sigma^{\textnormal{bot}}$};
		\fill[black] (0.1,0.63) circle (.075);
		\draw[gray,->] (-2.75,-2) --++(.5,-.2) node[right,yshift=-.1cm]{\small$y_1$};
		\draw[gray,->] (-2.75,-2) --++(.5,.2) node[right,yshift=.05cm]{\small$y_2$};
		\draw[gray,->] (-2.75,-2) --++(0,.5) node[above]{\small$x$};
	\end{tikzpicture}

	\caption{Left: schematic representation of the analysis of a fast-slow system near a cusp point (black dot located at the origin). Qualitatively speaking it can be shown that there exists an attracting manifold $\cS_\ve$ constructed from continuing $\cS_{\ve,-}$ beyond the cusp point (also black dot located at the origin). Right: $\Sigma^{\textnormal{bot}}$ is a small rectangle below the cusp point. It can be shown that: a) the image of  $\Sigma^{\textnormal{bot}}$ under the flow of \eqref{eq:cusp3} is an exponentially small strip about $\cS_\ve\cap\Sigma^\EX$ and b) an exponentially small interval $I_{\textnormal{bot}}\in\Sigma^{\textnormal{bot}}$ (shown in green) is mapped to an algebraically wide interval of order $\cO(\ve)$ in $\Sigma^\EX$. The main idea behind this observation is that we can choose trajectories that are exponentially close but that are separated by a repelling slow manifold. }
	\label{fig:cusp1}
\end{figure}
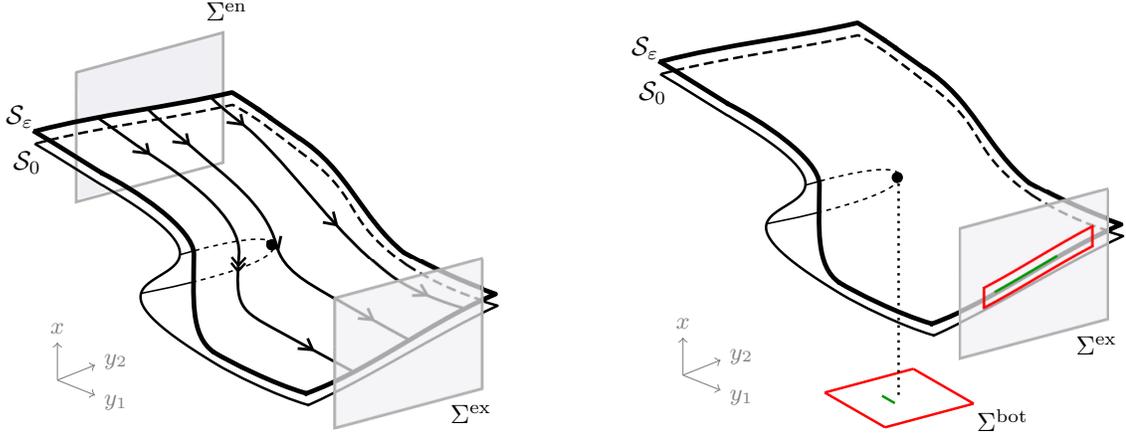

There are some refinements to the results of Theorem \ref{thm:cusp1} in the case where the higher order terms of \eqref{eq:cusp3} are of order $\cO(\ve)$. In particular, \cite{jardon2016analysis} provides more accurate information on the transition map $\Pi:\Sigma^{\EN}\to\Sigma^{\EX}$. The analysis performed in \cite{jardon2016analysis} uses a heavy amount of normal form theory, in particular results for quasihomogeneous vector fields \cite{lombardi2010normal}, an appropriate normal form for fast-slow systems \cite{jardon2015formal}, as well as for semi-hyperbolic singularities \cite{takens1971partially,bonckaert1997conjugacy,bonckaert1996partially}, and some other classical normal form techniques \cite{chen1963equivalence,sternberg1958structure}.\bigskip

\item[Bogdanov-Takens] Fast-slow systems related to a standard Bogdanov-Takens bifurcation have been studied using the blow-up method in \cite{de2011slow,chiba2011periodic}. Although the approaches are similar, the concerns and results are different, and here we sketch the most important ones. Let us start from \cite{de2011slow} where a planar fast-slow system of the form
\begin{equation}\label{eq:bt1}
	\begin{split}
			x' &= y,\\
			y' &= -xy + \ve(b_0+b_1x+x^2+x^3G(x,\omega))+y^2 H(x,y,\omega),
	\end{split}
\end{equation}
where $\omega$ denotes parameters in a compact subset of an Euclidean space (any finite dimension is allowed) and $b=(b_0,b_1)\sim 0\in\R^2$, are considered. It is important to mention that \eqref{eq:bt1} is a normal form, implying that the results of \cite{de2011slow} apply to a much wider class of planar fast-slow systems which after changes of coordinates can be brought to the form \eqref{eq:bt1}. Note that for $\ve>0$, \eqref{eq:bt1} indeed corresponds to (a versal unfolding of) the normal form of a Bogdanov-Takens bifurcation \cite{bogdanov1975versal} and \cite[Theorem 6.2]{takens2001forced}. The main result of \cite{de2011slow} consists on a full description of the bifurcation \emph{independent of $\ve$}, in contrast with \cite{bogdanov1975versal,takens2001forced} for which ``the validity of the results shrinks" as $\ve\to0$. The blow-up analysis is now considerably more involved than in the fold case, and is mostly performed in the family (or central) chart.

In \cite{chiba2011periodic} a one-parameter family of fast-slow systems presenting a Bogdanov-Takens point is studied, but in this case in a system with two fast and one slow variables of the form
\begin{equation}\label{eq:BT2}
	\begin{split}
		x_1' &= f_1(x_1,x_2,y,\ve,\delta),\\
		x_2' &= f_2(x_1,x_2,y,\ve,\delta),\\
		y'   &= \ve g(x_1,x_2,y,\ve,\delta).
	\end{split}
\end{equation}
Accordingly, the critical manifold is one-dimensional and is supposed to be  ``S-shaped", similar to the one appearing in the van der Pol oscillator. As such, the critical manifold has two fold points. The main assumption is that such fold points coincide with a Bogdanov-Takens bifurcation point of the layer equation. Furthermore, it is assumed that the stable part of the critical manifold correspond to stable foci of the layer equation, while the unstable part of the critical manifold corresponds to saddle nodes, see Figure \ref{fig:BT2} for a schematic representation of the critical manifold and the flow of the layer equation.

\begin{figure}[htbp]\centering
\begin{tikzpicture}
	\node at (0,0) {\includegraphics[scale=1.25]{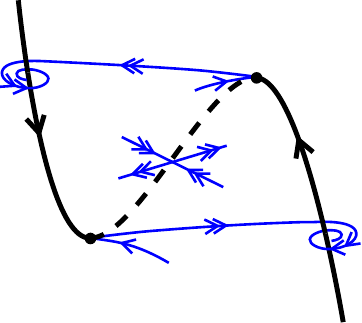}};
	\node at (-2,2.25) {$\cS_{0,+}^\txta$};
	\node at ( 2.1,-2.35) {$\cS_{0,-}^\txta$};
\end{tikzpicture}

	\caption{Schematic of the critical manifold of \eqref{eq:BT2} considered in \cite{chiba2011periodic}, where the critical manifold is a $1$-dimensional $S$-shaped curve. The parameter $\delta$, which controls the strength of attraction of $\cS_{0,\pm}^\txta$, plays a crucial role: for some pairs $(\ve,\delta)$ \eqref{eq:BT2} has a stable periodic orbit, while for other values of $(\ve,\delta)$ \eqref{eq:BT2} exhibits chaotic motion. }
	\label{fig:BT2}
\end{figure}

The main results of \cite{chiba2011periodic} prove and describe relaxation oscillations~\cite[Theorem 2]{chiba2011periodic} and chaotic motion~\cite[Theorem 3]{chiba2011periodic} of \eqref{eq:BT2} according to the parameters $(\delta,\ve)$. The parameter $\delta$ is used to control the strength of stability of the stable part of the critical manifold. The core analysis heavily depends on studying a first Painlev\'e equation (see \cite[Chapter 32]{olver2010nist} and \cite[Chapter 14.4]{ince1956ordinary}) in the blow-up space.\bigskip

\end{description}

The above examples are particular instances of fast-slow dynamical systems, where singularities beyond the classical fold appear. As it is evident, the analysis of fast-slow systems becomes more complicated, for example, when the number of fast variables increases and/or the degeneracy of the non-hyperbolic points also increases. There are several other situations that fall out of the scope of folded or Hopf-type singularities, see Sections \ref{sec:future} and \ref{sec:applications} for further references.

\subsection{Further directions}\label{sec:future}

In this survey we have reviewed research related to the foundations of the blow-up method as an important tool to analyze fast-slow systems. Naturally, there are some particular topics that have not been covered in this survey, and moreover, there are new emerging and exciting research directions, \ch{including works that consider problems beyond folded singularities} , where the blow-up method is being employed. The following are just some examples of works concerning and/or related to fast-slow dynamics with non-hyperbolic singularities, where the blow-up method has proved useful already:

\begin{description}[leftmargin=*]

	\item[Unbounded critical manifolds] One usually considers loss of normally due to tangency of the critical manifold and the fast foliation along sets of lower dimension than the critical manifold itself. Clear examples are the fold points and lines reviewed in Sections \ref{sec:genericfold}, \ref{sec:planarcanards}, \ref{sec:foldedsing} and cusp points sketched in Section \ref{sec:beyondfold}. However, it may occur that somewhere in the phase-space, the critical manifold aligns with the fast-foliation. These cases have been studied in \cite{GucwaSzmolyan,kuehn2014normal,kristiansen2017blowup}.

	\item[Codimension $0$ critical manifolds] In all examples provided so far, the critical manifold is at least of codimension $1$. \ch{In other words, the critical manifold is always of lower dimension than the phase-space.} For example, for the generic fold in Section \ref{sec:genericfold} the critical manifold is of codimension $1$ \ch{($\dim(\cC_0)=1$ and phase-space $\R^2$)}, while in the Bogdanov-Takens case of \cite{chiba2011periodic}, the critical manifold is of codimension $2$ \ch{($\dim(\cC_0)=1$ and phase-space $\R^3$)}. There are cases, however, where the critical manifold may have the same dimension as the phase-space. A trivial but non-generic example is given by the planar fast-slow system
\begin{equation}
	\begin{split}
		x' &= 0\\
		y' &= \ve,
	\end{split}
\end{equation}
where it is evident that the critical manifold is the the whole plane. A more elaborate example is given by
\begin{equation}
	\begin{split}
		x' &= f(x,y,\ve)\\
		y' &= \ve g(x,y,\ve),
	\end{split}
\end{equation}
where
\begin{equation}
	f(x,y,\ve)=\begin{cases}
		h(y)\exp\left( -\frac{1}{x} \right)+ \cO(\ve), & x>0\\
		 \cO(\ve), & x\leq0,
	\end{cases}
\end{equation}
where $h(y)$ can be any smooth nonzero function. The function $f(x,y,\ve)$ is smooth and the critical manifold is given by
\begin{equation}
	\mathcal C_0 = \left\{ (x,y)\in\mathbb R^2\, | \, x\leq 0 \right\},
\end{equation}
which is two dimensional. Similarly, many other examples can be created.

	In a more relevant context, this situation of codimension $0$ critical manifolds has recently been considered in \cite{kuehn2018duck} motivated by hysteresis operators~\cite{kuehn2017generalized}.

	\item[Piecewise-smooth systems] Canard theory (see Sections \ref{sec:planarcanards} and \ref{sec:foldedsing}) has been extended to piecewise-linear systems~\cite{desroches2016canards,desroches2018piecewise} and in a similar context the blow-up method has been used to analyze piece-wise smooth (fast-slow) systems \cite{kristiansen2015use,kristiansen2015regularizations,de2018non,jeffrey2016hidden}.

	\item[Torus canards] In all the fast-slow systems considered so far, the critical manifold is made of critical points of the fast dynamics. However, in high dimensional problems, the dynamics of the layer equation may converge to, for example, limit cycles instead of equilibrium points. In this context, canard theory has been extended to cases where the fast dynamics has limit cycles as limiting sets, leading to the so-called torus canards~\cite{burke2016canards,vo2017generic}, see also \cite{ju2018bottom}.

	\item[Discrete time systems] The blow-up method has been recently used to study discrete fast-slow systems with a transcritical singularity \cite{engel2018discretized} directly without using continuous approximations. This represents additional steps to extending the blow-up to the context of maps with singularities; see also~\cite{NippStoffer2,NippStofferSzmolyan}.
\end{description}

\subsection{Applications}\label{sec:applications}

The main part of this document has been dedicated to reviewing important theoretical progress on the understanding of fast-slow dynamics using the blow-up method. As such, we have emphasized the usefulness of the method to, for example, explain the canard explosion phenomenon. However, the blow-up method has also been used in interesting particular applications. Here we briefly summarize some examples to illustrate the breadth of different areas, where the blow-up technique can be useful.

\begin{description}[leftmargin=*]
	\item[In neuroscience] The progress in the theory of fast-slow systems, particularly the understanding of canards, and the related mixed-mode oscillations~\cite{desroches2012mixed}, has had a great impact in neuroscience. As an example, the Fitzugh-Nagumo, the Hodgkin-Huxley, and other related models have been studied from several mathematical perspectives, especially from a multi-time scale approach \cite{Guckenheimer2009,Guckenheimer2010,Desroches2013,Rubin2007}. A recent instance where the blow-up method is instrumental is \cite{Carter2018unpeeling}, where the transition between two distinct oscillatory patters in the Fitzugh-Nagumo system is precisely described.

 \item[In systems biology] The analysis of a model for the division cycle of an embryonic cell, known as the \emph{Goldbeter model}, is studied in \cite{kosiuk2016geometric}. This is the first rigorous analysis of a system based on Michaelis-Menten equations where small Michaelis-Menten constants are considered and play a central role. In its original form, the model under study has no evident time scale separation. After a transformation, an auxiliary and equivalent nonstandard fast-slow model is obtained. By using tools of geometric singular perturbation theory the authors analyze the system using the property that the auxiliary model has three time scales. On a related context, a problem exhibiting relaxation oscillations and depending on two singular parameters is considered in~\cite{GucwaSzmolyan2}.

\item[In chemistry]	The \emph{Olsen model} for the peroxidase–oxidase reaction is considered in \cite{kuehn2015multiscale}. There the model is $4$-dimensional. Besides this evident obstacle, the Olsen model is not in standard form \eqref{eq:sfs1}, has several small parameters, and has non-fold singularities. Even more interestingly, the Olsen model has no defined fast and slow variables, but presents three distinct regimes in which there are $1$, $2$, and $3$ fast variables. The blow-up method is essential in this case to understand the dynamics in all the regimes and leads to the description of non-classical relaxation oscillations of the Olsen model.

	\item[In engineering] In \cite{Bossolini2017}, a model describing an earthquake fault is considered. The model being studied is a Hamiltonian system, and possesses an \emph{unbounded critical manifold}, see further information above. The blow-up method is used to show the existence of limit cycles associated to a degenerate Hopf bifurcation. Another recent application is found in \cite{iuorio2018}, where a Micro-Electro Mechanical System (MEMS) undergoing a singular effect known as ``touchdown'' is analyzed. In this case, the blow-up method allows a precise analysis of the steady states of a regularized model, and the description of the dynamics as a couple of small parameters tend to zero.

	\item[In control theory] Singular perturbation problems are classical in control theory. The case for which the critical manifold is normally hyperbolic is in fact well understood \cite{kokotovic1999singular}. The case where the critical manifold looses normal hyperbolicity had not received much attention. Recently, the blow-up method has allowed the design of controllers that stabilize non-hyperbolic points (of any degeneracy) of fast-slow systems with one fast variable \cite{jardon2018improving,jardon2019stabilization}, and the model order reduction of fast-slow systems with hidden/degenerate normally hyperbolic critical manifolds \cite{jardon2017model}.

\end{description}

\section{Summary and Outlook}\label{sec:summary}

This survey has been devoted to provide a concise recollection of the most influential progress in the blow-up theory of fast-slow systems. As such, we have principally covered folded singularities and in particular how the blow-up method allowed a geometric description of canards, canard cycles, and canard explosion. We have also outlined on the role of the blow-up method when analyzing other singularities beyond the fold, such as in the case of the Hopf, transcritical, pitchfork, cusp, and Bogdanov-Takens points. Moreover, we have briefly reviewed some recent theoretical and applied research in fast-slow systems, where the blow-up method is also one of the main mathematical tools for analysis.

It is evident that the blow-up method enables us to study and understand complicated fast-slow dynamics around non-hyperbolic points. However, it should also be clear that the difficulties do not end once we apply the method. In several cases, the dynamics within the charts are still quite intricate and, for example, one encounters special equations for which their asymptotics are important (e.g. Riccati, Weber, and Panilev\'e). Furthermore, one must trace invariant manifolds, compute their distance (Melknikov integral), control the obstacle of resonances, or deal with complex and/or complicated path integrals, just to mention a few of the difficulties.

We predict that for many of the theoretical and applied studies of fast-slow systems to come, the blow-up method will likely continue to be one of the techniques of choice when studying such dynamical systems. Furthermore, due to the everyday increase of complexity of systems and their models, the blow-up method will have to be extended and or adapted to such needs. We envision that, for instance, geometric studies of high dimensional systems and/or high codimension singularities are going to increase in relevance. Finally, although the blow-up method has been a great success for deterministic ODE systems, it is at this point a daunting challenge to extend it to a wider range of dynamical systems including partial and stochastic or even stochastic-partial differential equations.

\begin{center}
	\vspace*{.25cm}
	\textbf{Acknowledgments}
\end{center}
The authors gratefully acknowledge the comments of the anonymous reviewer that helped to improve the content and presentation of this article. HJK is supported by a Technical University Foundation Fellowship and by an Alexander von Humboldt Research Fellowship. CK would like to thank the VolkswagenStiftung for support via a Lichtenberg Professorship.

\newcommand{\etalchar}[1]{$^{#1}$}
\providecommand{\bysame}{\leavevmode\hbox to3em{\hrulefill}\thinspace}
\providecommand{\MR}{\relax\ifhmode\unskip\space\fi MR }
\providecommand{\MRhref}[2]{%
  \href{http://www.ams.org/mathscinet-getitem?mr=#1}{#2}
}
\providecommand{\href}[2]{#2}

\end{document}